\documentclass[10pt,leqno]{article}
\usepackage{amsmath,amsthm}
\usepackage{amssymb}
\usepackage[all]{xypic}
\usepackage{amsbsy}

\swapnumbers
\newtheorem{thm}{Theorem}[section]
\newtheorem{lem}[thm]{Lemma}
\newtheorem{prop}[thm]{Proposition}
\newtheorem{cor}[thm]{Corollary}
\newtheorem{defn}[thm]{Definition}
\theoremstyle{definition}
\newtheorem{rem}[thm]{Remark}
\newtheorem{exam}[thm]{Example}

\numberwithin{equation}{section}
\numberwithin{figure}{section}

\def\Aut{\operatorname{Aut}}

\def\ZZ{{{\mathbb{Z}}}}

\def\QQ{{{\mathbb{Q}}}}

\def\cE{{\cal E}}
\def\cF{{\cal F}}

\def\cM{{\cal M}}

\def\cU{{\cal U}}

\def\La{{\Lambda}}

\def\colim{\operatorname{colim}}
\def\holim{\operatorname{holim}}

\def\End{\mathrm{End}}
\def\Ext{{{\mathrm{Ext}}}}
\def\Ga{\Gamma}

\def\Hom{\mathrm{Hom}}

\def\Mod{\mathbf{Mod}}
\def\map{\operatorname{map}}

\def\longr{{\longrightarrow}}
\def\defeq{\overset{\mathrm{def}}=}
\def\mm{{{\mathfrak{m}}}}

\def\Spec{{{\mathrm{Spec}}}}
\def\GG{{{\mathbb{G}}}}
\def\FF{{{\mathbb{F}}}}

\def\End{{{\mathrm{End}}}}

\newcommand{\brackets}[4]{{{\left\{ \begin{array}{ll}
				{{#1}}&{{#2}}\\
				\\
				{{#3}}&{{#4}}\end{array}\right.}}}

\def\fgl{{\mathbf{fgl}}}

\def\cmpl{{\wedge}}

\def\rege{{{\epsilon}}}

\date{}

\def\CP{{{\mathbb{C}\mathrm{P}}}}

\def\Ell{{{e\ell\ell}}}
\def\\tmf{{{\mathbf{\tmf}}}}

\def\tmf{{{\mathbf{tmf}}}}

\def\Gal{{{\mathrm{Gal}}}}
\def\comod{{{\mathbf{Comod}}}}
\def\fglp{{\mathbf{fglp}}}

\begin{document}

\title{The Adams-Novikov Spectral Sequence\\ and\\
the Homotopy Groups of Spheres}
\author{Paul Goerss\footnote{The
author were partially supported by the National Science
Foundation (USA).}}
\maketitle
\begin{abstract}These are notes for a five lecture series intended to
uncover large-scale phenomena in the homotopy groups of
spheres using the Adams-Novikov Spectral Sequence. The
lectures were given in Strasbourg, May 7--11, 2007.\end{abstract}
\bigskip
\centerline{\today}
\bigskip

\def\msing{{{\cM_\Ell^{\mathrm{sing}}}}}
\def\mshat{{{\widehat{\cM}_\Ell^{\mathrm{sing}}}}}
\def\msmooth{{{\cM_\Ell^{\mathrm{sm}}}}}
\def\mweier{{{\cM_{\mathrm{Weier}}}}}
\def\msmi{{{\cM_{\mathrm{SMI}}}}}

%{\bf Warning!} These notes are in extremely preliminary form and are
%really intended only as a supplement to the actual lectures. They
%will be updated continuously between now and the end of the lecture
%series; therefore, any user should pay attention to the date of the version
%they have -- and all experts should not get too exercised about the
%sources, as they too will be updated. Of course, if you spot errors,
%let me know.
%\bigskip

\tableofcontents
\bigskip

\bigskip

{\bf A note on sources:} I have put some references at
the end of these notes, but they are nowhere near exhaustive.
They do not, for example, capture the role of Jack Morava in
developing this vision for stable homotopy theory. Nor somehow,
have I been able to find a good way to record the overarching
influence of Mike Hopkins on this area since the 1980s. And, although,
I've mentioned his name a number of times in this text, I also
seem to have short-changed Mark Mahowald -- who, more than
anyone else, has a real and organic feel for the homotopy groups
of spheres. I also haven't been
very systematic about where to find certain topics. If I seem a bit
short on references, you can be sure I learned it from the absolutely
essential reference book by Doug Ravenel \cite{Rav} -- ``The Green Book'',
which is not green in its current edition. This book contains
many more references. Also essential is the classic
paper of Miller, Ravenel, and Wilson \cite{MRW} and the papers those
authors wrote to make that paper go -- many of which I haven't cited
here. Another source is Steve Wilson's sampler \cite{primer}, which also has some
good jokes. 

\section{The Adams spectral sequence}

Let $X$ be a spectrum and let $H^\ast (X) = H^\ast (X,\FF_2)$
be its mod $2$ cohomology. Then $H^\ast X$ is naturally
a left module over the mod $2$ Steenrod algebra $A$.

We will assume for the moment
that all spectra are {\it bounded below} and of {\it finite type};
that is $\pi_iX = 0$ if $i$ is sufficiently negative and $H^nX$ is
a finite vector space for all $n$. These are simply pedagogical
assumptions and will be removed later. 

The {\it Hurewicz map} is the natural homomorphism
\begin{align*}
[X,Y] &\longr \Hom_A(H^\ast Y,H^\ast X)\\
f& \mapsto f^\ast = H^\ast(f).
\end{align*}
and we would like to regard the target as relatively computable.
However, the map is hardly ever an isomorphism.

There is, however, a distinguished class of spectra for which it {\it is}
 an isomorphism.
Suppose $V$ is a graded $\FF_2$-vector space, finite in each
degree and bounded below; then by
Brown representability, there is a spectrum 
$KV$ -- a generalized mod $2$ {\it Eilenberg-MacLane spectrum} --
and a natural isomorphism
$$
[X,KV] \cong \Hom_{\FF_2}(H_\ast X,V).
$$
Such a spectrum has the following properties. 
\begin{enumerate}

\item There is a natural isomorphism $\pi_\ast KV \cong V$ (set
$X = S^n$ for various $n$);

\item $H^\ast KV$ is a projective $A$-module -- indeed, any splitting
of the projection $H^\ast KV \to V^\ast$ defines
an isomorphism $A \otimes V^\ast \to H^\ast KV$.
\end{enumerate}
Then for all $X$ (satisfying our assumptions) we have
$$
[X,KV] \cong \Hom_{\FF_2}(V^\ast,H^\ast X) \cong 
\Hom_A(H^\ast KV,X).
$$
Let me isolate this as a paradigm for later constructions:

\begin{lem}\label{basic}There is a natural homomorphism
$$
h:[X,Y] \longr \Hom_A(H^\ast X,H^\ast Y)
$$
which is an isomorphism when $Y$ is a generalized mod $2$
Eilenberg-MacLane spectrum. Furthermore, for such $Y$, $H^\ast Y$
is projective as an $A$-module.
\end{lem}

The idea behind the Adams Spectral Sequence (hereinafter know
as the clASS) and all its variants is that we resolve a general $Y$ by
Eilenberg-MacLane
spectra to build a spectral sequence for computing $[X,Y]$ with
edge homomorphism $h$.

\begin{defn}\label{adams-res}Let $Y$ be a spectrum. Then
an Adams resolution for $Y$ is a sequence of spectra
$$
\xymatrix@C=15pt{
Y \rto^-d& K^0 \rto^-d &K^1 \rto^-d& K^2 \to \cdots
}
$$
where each $K^s$ is an mod $2$ Eilenberg-MacLane spectrum and for all 
mod $2$ Eilenberg-MacLane spectra, the chain complex of
vector spaces
$$
\cdots \to [K^2,K] \to [K^1,K] \to [K^0,K] \to [Y,K] \to 0
$$
is exact.
\end{defn}

Such exist and are unique up to a notion of chain equivalence which I 
leave  you to formulate. It also follows  -- by setting $K = \Sigma^n H\ZZ_2$
for various $n$ --  that
$$
\cdots \to H^\ast K^1 \to H^\ast K^0 \to H^\ast Y \to 0
$$
is a projective resolution of $H^\ast Y$ as an $A$-module.

\begin{rem}[{\bf The Adams Tower}]\label{adams-tower}
From any Adams resolution of $Y$ we can build a tower of spectra
under $Y$ 
$$
\cdots \longr Y_2 \longr Y_2 \longr Y_0
$$
with the properties that
\begin{enumerate}

\item $Y_0 = K^0$ and $Y \to Y_0$ is the given map $d$;

\item there are homotopy pull-back squares for $s > 0$ 
$$
\xymatrix{
Y_s \rto \dto &\ast\dto\\
Y_{s-1} \rto & \Sigma^{1-s}K^s;
}
$$

\item the induced maps
$$
\Sigma^{-s}K^s \to Y_{s+1} \to \Sigma^{-s}K^{s+1}
$$
are the given maps $d$; and

\item the induced map $\Sigma^{-s} K^s \to Y_s$ induces
a short exact sequence
$$
0 \to H^\ast Y \to H^\ast Y_s \to \Sigma^{-s} M_s \to 0
$$
where 
$$
M_s = \mathrm{Ker}\{d^\ast:H^\ast K^s \to H^\ast K^{s-1}\} =
\mathrm{Im}{d^\ast: H^\ast K^{s+1} \to H^\ast K^s}.
$$
The short exact sequence is split by the map $Y \to Y_s$.
\end{enumerate}
\end{rem}

From the tower (under $Y$ remember) we get an induced
map
$$
\colim H^\ast Y_s \to H^\ast \lim Y_s \to H^\ast Y
$$
and the composite is an isomorphism. Under my hypotheses
of bounded below and finite-type, the first map is an isomorphism
as well; from this and from the fact that homotopy inverse limits
of $H_\ast$-local spectra are local, it follows that $Y \to \lim Y_s$
is the $H_\ast$-localization of $Y$, which I will write $Y^\cmpl_2$,
for if the homotopy groups of $Y$ are finitely generated, the homotopy
groups of the localization are gotten by completion at $2$ from
the homotopy groups of $Y$.

We now get a spectral sequence built from the exact couple
$$
\xymatrix{
[\Sigma^{t-s}X,Y_{s+1}] \rto & [\Sigma^{t-s}X,Y_s] \rto \ar@{-->}[dl]
& [\Sigma^{t-s},Y_{s-1}]  \ar@{-->}[dl]\\
[\Sigma^{t-s}X,\Sigma^{-s-1}K^{s+1}]\ar[u] &[\Sigma^{t-s}X,\Sigma^{-s}K^s]\ar[u]
}
$$
where the dotted arrows are of degree $-1$. This spectral
sequence has $E_1$ term
\begin{align*}
E_1^{s,t} &\cong [\Sigma^{t-s}X,\Sigma^{-s}K^s]\\
&\cong [\Sigma^{t}X, K^s]\\
&\cong \Hom_A(H^\ast K^s,\Sigma^tH^\ast X)
\end{align*}
and from this and Remark \ref{adams-tower}.3, we have that
\begin{align*}
E_2^{s,t}&\cong H^s\Hom_A(H^\ast K^\bullet,\Sigma^t H^\ast X)\\
&\cong \Ext_A^s(H^\ast Y,\Sigma^tH^\ast X).
\end{align*}
It abuts to $[\Sigma^{t-s}X,\lim Y_s]$; thus we write the clASS
as
$$
\Ext_A^s(H^\ast Y,\Sigma^tX) \Longrightarrow [\Sigma^{t-s}X,Y^\cmpl_2].
$$
The differentials go
$$
d_r:E_r^{s,t} \to E_r^{s+r,t+r-1}.
$$
This is a left half-plane spectral sequence. For large $r$, we have
$$
E_\infty^{s,t} \subseteq E_r^{s,t}
$$
but it can happen that there is no finite $r$ for which this is
an equality. Thus we need to distinguish between {\it abuts}
and {\it converges}.

Write $e_\infty^{s,t}$ as the quotient of
$$
\mathrm{Ker}\{[\Sigma^{t-s}X,Y^\cmpl_2] \to [\Sigma^{t-s}X,Y_{s-1}]\}
$$
by
$$
\mathrm{Ker}\{[\Sigma^{t-s}X,Y^\cmpl_2] \to [\Sigma^{t-s}X,Y_{s}]\}.
$$
Then, all the word ``abuts'' signifies is that there is a map
$$
e_\infty^{s,t} \longr E_\infty^{s,t}
$$
where the target is the usual $E_\infty$ term of the spectral sequence.

\begin{defn}\label{converges}
The spectral sequence {\bf converges} if
\begin{enumerate}

\item this map is an isomorphism; and

\item $\lim^1 [\Sigma^nX,Y_s] = 0$ for all $n$.
\end{enumerate}
\end{defn}

We get convergence if, for example, we know that 
for all $(s,t)$ there is an $r$ so that
$E_\infty^{s,t} = E_r^{s,t}$. In particular, this will happen (although
this is not obvious) if $H_\ast X$ is {\it finite} as a graded vector
space. See Theorem \ref{vanishing} below. A thorough discussion
of convergence and the connections  with localization is
contained in Bousfield's paper \cite{stlocal}. A source for
the Adams spectral sequence and the Adams-Novikov Spectral
Sequence is \cite{blue}, although the point of view I'm adopting
here is from the first section of \cite{moore}.

\section{Classical calculations}

 In any of the Adams spectral sequences, there are two computational
steps: (1) the algebraic problem of calculating the $E_2$-term; and
(2) the geometric problem of resolving the differentials. Both
problems are significant and neither has been done completely.
There has been, at various times, optimism that resolving the first
would lead to a resolution of the second, but the following
principle, first named by Doug Ravenel, seems to hold
universally:
\bigskip

{\bf The Mahowald Uncertainty Principle:} Any spectral sequence
converging to the homotopy groups of spheres with an $E_2$-term that
can be named using homological algebra will be infinitely far from
the actual answer.
\bigskip

Nonetheless, we soldier on. For the clASS we have at least three
tools for computing the $E_2$-term of
$$
\Ext^s_A(H^\ast Y,\Sigma^t\FF_2) \Longrightarrow \pi_{t-s}Y^\cmpl_2.
$$
They are
\begin{enumerate}

\item minimal resolutions, which can be fully automated, at least at
the prime $2$ (see \cite{bruner});

\item the May spectral sequence (still the best method for by-hand
computations; see \cite{tangora} and \cite{Rav} \S 3.2); and

\item the $\Lambda$-algebra, which picks up the interplay with
unstable phenomena (see \cite{6A});
\end{enumerate}
I won't say much about them here, but all of them can be a
lot of fun.

We now concentrate on calculating the clASS
$$
E_2^{s,t}(Y)=\Ext_A(H^\ast Y,\Sigma^t\FF_2) \Longrightarrow
\pi_{t-s}Y
$$
where $Y$ is a homotopy associative and homotopy commutative
$H_\ast$-local
{\it ring} spectrum. For such $Y$ the the homotopy groups of $Y$ are
a graded commutative ring and the Adams spectral sequence
becomes a spectral sequence of appropriate graded commutative
rings. In particular, differentials are derivations. There are also
Massey products and Toda brackets. 

Let me say a little more about these last. Suppose
$a,b,c \in \pi_\ast Y$ have degree $m$, $n$, and $p$ respectively
and suppose $ab = 0 = bc$. Then, by choosing null-homotopies
we get a map
$$
\langle a,b,c \rangle:S^{m+n+p+1} =
C^{m+n+1} \wedge S^p \coprod_{S^m\wedge S^n \wedge S^p}
S^m \wedge C^{n+p+1} \to Y
$$
where $C^{n+1}$ is the cone of $S^n$. This bracket is
well-defined up to the choices of null-homotopy; that is, well-defined
up to the indeterminacy
$$
a[\pi_{n+p+1}Y] + [\pi_{m+n+1}Y]c.
$$
There is a similar construction algebraically; that is, given
$$
x,y,z \in \Ext_A^\ast(H^\ast Y,\Sigma^\ast \FF_2)
$$
with $xy = 0$ and $yz=0$, there is a Massey product
$$
\langle x,y,z \rangle  \in \Ext_A^\ast(H^\ast Y,\Sigma^\ast \FF_2)
$$
well-defined up to indeterminacy. Furthermore, if the Massey product
detects the Toda bracket via the spectral sequence.

\begin{rem}[{\bf The cobar complex}]\label{bar-complex}
The Massey product is computed
when the $\Ext$ groups are given as the cohomology of
a differential graded algebra. We now describe the usual
example of such a dga.

The
forgetful functor from $A$-modules to graded vector spaces
has a left adjoint $V \mapsto A \otimes V$. The resulting
cotriple resolution of an $A$-module $M$ is the {\it bar} complex
$B_\bullet (M)\to M$. This is an augmented simplicial $A$-module with
$$
B_n(A) = A^{\otimes (n+1)} \otimes M
$$
and face operators given by
$$
d_i(a_0 \otimes \ldots a_n \otimes x) = \brackets{a_0 \otimes 
\ldots \otimes a_ia_{i+1} \otimes \ldots a_n \otimes x,}{i < n;}
{a_0 \otimes \ldots \otimes a_{n-1} \otimes a_nx,}{i=n.}
$$
The map
$$
a_0 \otimes \ldots a_n \otimes x \to 1 \otimes a_0 \otimes \ldots a_n \otimes x
$$
gives this complex a contraction as a simplicial vector space; hence we
have a natural projective resolution. Applying $\Hom_A(-,\Sigma^\ast\FF_2)$
to the bar complex yields the cobar complex $C^\bullet(M_\ast)$
where 
$$
C^n(M_\ast) = A_\ast^{\otimes n} \otimes M_\ast
$$
where I have written $A_\ast$ and $M_\ast$ for the duals. The
coface maps now read
$$
d^i(\alpha_1 \otimes \ldots \alpha_n \otimes y)
={{{\left\{ \begin{array}{ll}
				{{1\otimes \alpha_1 \otimes \ldots \alpha_n \otimes y,}}&
				{{i=0;}}\\
				\\
				{{\alpha_1 \otimes \ldots \otimes \Delta(\alpha_i)
				\otimes \ldots \otimes \alpha_n\otimes y,}}&{{0 <i<n+1;}}\\
				\\
				{{\alpha_1 \otimes \ldots \otimes \alpha_n 
				\otimes \psi(y).}}&{{i=n+1.}}\end{array}\right.}}}
$$
Here $\Delta:A_\ast \to A_\ast \otimes A_\ast$ and 
$\psi:M_\ast \to A_\ast \otimes M_\ast$ are the diagonal and
comultiplication respectively. If $M_\ast =H_\ast Y$,
where $Y$ is a ring spectrum, then $C^\bullet (H_\ast Y)$ is
a cosimplicial $\FF_2$-algebra and
$$
\pi^\ast C^\bullet(H_\ast Y) = \Ext^\ast_A(H^\ast Y,\Sigma^\ast \FF_2).
$$
For computational purposes, we often take the normalized
bar construction $\bar{C}^\bullet(M_\ast)$ where
$$
\bar{C}^n(M_\ast) \cong \bar{A}_\ast^{\otimes n} \otimes M_\ast
$$
where $\bar{A}_\ast \subseteq A_\ast$ is the augmentation
ideal. If $M=H_\ast Y$, then $\bar{C}^\bullet(H_\ast Y)$
is a differential graded algebra and a good place to compute
Massey products.
\end{rem}

Some beginning computational results are these.

\begin{thm}\label{0-and-1-lines}1.) $\Hom_A(\FF_2,\FF_2) = \FF_2$ generated
by the identity. If $t \ne 0$, then $\Hom_A(\FF_2,\Sigma^t\FF_2) = 0$.

2.) $\Ext_A^1(\FF_2,\Sigma^t\FF_2) \cong \FF_2$ if $t=2^i$ for some $i$;
otherwise this group is zero.
\end{thm}

\begin{proof}The short exact sequence of $A$-modules
$$
0 \to \bar{A} \to A \to \FF_2 \to 0
$$
and the resulting long exact sequence in $\Ext$ shows that
\begin{align*}
\Ext_A^1(\FF_2,\Sigma^\ast\FF_2) &\cong \Hom_A(\bar{A},\Sigma^\ast\FF_2)\\
&\cong \Hom_{\FF_2}(\FF_2\otimes_A \bar{A},\Sigma^\ast \FF_2)\\
&\cong \Hom_{\FF_2}(\bar{A}/\bar{A}^2,\Sigma^\ast \FF_2)
\end{align*}
and the Adem relations show that $\bar{A}/\bar{A}^2$ is zero
except in degrees $2^i$ where it is generated by the residue class
of $\mathrm{Sq}^{2^i}$.
\end{proof}

The non-zero class is $\Ext_A^1(\FF_2,\Sigma^{2^i}\FF_2)$ is called
$h_i$.

The following result shows that the clASS for the sphere actually converges.
See Definition \ref{converges} and the remark following. Part (3) is usually phrased as saying there is a vanishing line of slope $2$.
The following result was originally due to Adams \cite{periodicity}.

\begin{thm}\label{vanishing}We have the following vanishing result
for the clASS:
\begin{enumerate}

\item $Ext^s_A(\FF_2,\Sigma^{t}\FF_2) = 0$ all $(s,t)$ with $t-s<0$.

\item $Ext^s_A(\FF_2,\Sigma^{s}\FF_2) = \FF_2$ generated by $h_0^s$;

\item $Ext^s_A(\FF_2,\Sigma^{t}\FF_2) = 0$ for all $(s,t)$ with
$$
0 < t-s < 2s + \rege
$$
where
$$
\rege ={{{\left\{ \begin{array}{ll}
				{{1,}}&
				{{s\equiv 0,1\ \mathrm{mod}\ (4);}}\\
				\\
				{{2,}}&{{s\equiv 2\ \mathrm{mod}\ (4);}}\\
				\\
				{{3,}}&{{s\equiv 3\ \mathrm{mod}\ (4);.}}\end{array}\right.}}}
$$
\end{enumerate}
\end{thm}

We have the following low dimensional calculations. The first
is from \cite{hopf}.

\begin{thm}\label{2-line-ord}The $2$-line $Ext^2_A(\FF_2,\Sigma^\ast\FF_2)$
is the graded vector space generated by $h_ih_j$ subject to the
relations $h_ih_j = h_jh_i$ and $h_ih_{i+1}=0$.
\end{thm}

The relation in this last result immediately allows for Massey products;
a fun exercise is to show 
$$
\langle h_0,h_1 h_0 \rangle = h_1^2 \qquad \mathrm{and}\qquad
\langle h_1,h_0, h_1 \rangle = h_0h_2.
$$
Deeper results include the following from \cite{wang}:

\begin{thm}\label{3-line}The $3$-line $Ext^3_A(\FF_2,\Sigma^\ast\FF_2)$
is the graded vector space generated $h_ih_jh_k$ subject to the
further relations
$$
h_ih_{i+2}^2 = 0\qquad \mathrm{and}\qquad  h_i^2h_{i+2} = h_{i+1}^2
$$
and the classes
$$
c_i =\langle h_{i+1},h_i,h_{i+2}^2\rangle
$$
\end{thm}

\begin{rem}\label{some-classes}Surprisingly few of the low dimensional
classes actually detect homotopy classes. The classes $h_i$, $0 \leq i \leq
3$ detect $\times 2:S^0 \to S^0$ and the stable Hopf maps
$\eta$, $\nu$, and $\sigma$ respectively. Otherwise we have,
from \cite{hopf}, that
$$
d_2 h_i = h_{i-1}h_0^2.
$$
The first of these, when $i=4$ is forced by the fact that $2\sigma^2=0$,
as it is twice the square of an odd dimensional class in a graded commutative
ring. The others can be deduced from Adams's original work on secondary
operations or by  a clever inductive argument, due to Wang \cite{wang},
which relies on knowledge of the $4$-line.

On the $2$-line, the classes $h_0h_2$, $h_0h_3$, $h_2h_4$,
$h_i^2$, $i \leq 5$, and $h_jh_1$, $j \geq 3$ survive and detect homotopy
classes. The infinite family detected by the classes $h_jh_1$ is
the family known as Mahowald's $\eta_j$-family \cite{etaj}; they provided
a counterexample to the ``Doomsday conjecture'' -- which posited that there
could only be finitely many non-bounding infinite cycles on any $s$-line of the
clASS. 

The only other classes on the $2$-line which may survive are the Kervaire
invariant classes $h_i^2$, $i \geq 6$. This remains one of the more
difficult problems in stable homotopy theory. If you want to
work on this question
make sure you are in constant contact with Mark Mahowald, Fred Cohen,
and/or Norihiko
Minami. Preferably ``and'', not ``or''. There are many standard errors 
and seductive blind alleys; an experienced guide is essential. See
\cite{minami} for many references. One
world view suggests that not only is $h_i^2$ a permanent cycle,
but it is also the stabilization of a class $x\in \pi_\ast S^{2^{i+1}-1}$
with the property that $2x$ is the Whitehead product. It is this last
problem we would really, really like to understand.\footnote{In my youth, I once
suggested to Fred Cohen that this problem wasn't so important any more.
He became visibly upset and suggested I didn't understand homotopy theory.
He was right -- even allowing for the self-importance of the new guy in town,
this was a remarkably misguided statement. I apologize, Fred.}

\end{rem}

The first large scale family of interesting classes detecting homotopy
elements is in fact near the vanishing line. We have that $h_0^4h_3=0$,
so if 
$$
x \in \Ext_A^s(\FF_2,\Sigma^t\FF_2)
$$
with $x$ so that $h_0^4x$ is above the vanishing line, we can form
the Massey product $\langle x,h_0^4,h_4\rangle$. The following
result is the first result of Adams periodicity \cite{periodicity}.

\begin{thm}\label{adams-1}The operator
$$
P(-) = \langle -,h_0^4,h_3\rangle:
\Ext_A^s(\FF_2,\Sigma^t\FF_2) \to
\Ext_A^{s+4}(\FF_2,\Sigma^{t+12}\FF_2)
$$
is an isomorphism whenever $(s+4,t+4)$ is above the vanishing line.
\end{thm}

There is actually a stronger result. For all $n > 1$, we have
$h_0^{2^{n}}h_{n+1} = 0$ and, again near the vanishing line, we
can form the operator
$$
\langle -,h_0^{2^n},h_{n+1}\rangle:
\Ext_A^s(\FF_2,\Sigma^t\FF_2) \to
\Ext_A^{s+2^n}(\FF_2,\Sigma^{t+3\cdot 2^n}\FF_2)
$$
and this tends to be an isomorphism where it is universally defined, which
is approximately above a line of slope $1/5$. Thus in this ``wedge'',
between the lines of slope $1/5$ and slope $1/2$, we have considerable
hold over the clASS. 

The operator $P(-)$ of Theorem \ref{adams-1} takes permanent cycles
to permanent cycles; however, there is no guarantee that they are not
boundaries -- indeed, being so high up there are lot of elements that
can support differentials below. In order to show they are not boundaries,
we need to use other homology theories.

\section{The Adams-Novikov Spectral Sequence}

In working with generalized (co-)homology theories, it  turns out
to be important to work with homology rather than cohomology. This
is even evident for ordinary homology -- cohomology is the dual of
homology (but not the other way around) and cohomology is actually,
thus, a graded topological vector space. This is the reason for the finite
type hypotheses of the previous sections. There is no simple hypothesis
to eliminate this problem
for generalized homology theories: $E^\ast E$ will almost always
be a topological $E^\ast$-module. The price we pay is that we must
work with comodules rather than modules and, while the two notions
are categorically dual, the notion of a comodule is much less intuitive
(at least to me) than that of a module.

Let $E_\ast$ be a homology theory. We make many assumptions about
$E$ as we go along. Here's the first two.
\begin{enumerate}

\item The  associated cohomology theory $E^\ast$ has graded commutative
cup products. Thus the representing spectrum $E$ is a homotopy associative
and commutative ring spectrum.

\item There are two morphisms, the left and right units
$$
\eta_L,\eta_R:E_\ast \longr E_\ast E
$$
interchanged by the conjugation $\chi:E_\ast E \to E_\ast E$.
We assume that for one (and hence both) of these maps,
$E_\ast E$ is a flat $E_\ast$-module.
\end{enumerate}

The second of these conditions eliminates some of the standard
homology theories: $H_\ast (-,\ZZ)$ and connective $K$-theory, for example.

Under these assumptions, there is a natural isomorphism of
left $E_\ast$-modules
$$
E_\ast E \otimes_{E_\ast} E_\ast X \mathop{\longr}^{\cong}
E_\ast (E \wedge X)
$$
sending $a \in E_m E$ and $x \in E_n X$ to the composition
$$
\xymatrix{
S^n \wedge S^m \rto^-{a \wedge x}& E \wedge E \wedge E \wedge X
\ar[rr]^-{E \wedge m \wedge X} && E \wedge E \wedge X
}
$$
Note that the right $E_\ast$-module structure on $E_\ast E$ gets
used in the tensor product. In particular we have
$$
E_\ast (E \wedge E) \cong E_\ast E \otimes_{E_\ast} E_\ast E
$$
and the map
$$
E \wedge \eta \wedge E: E \wedge S^0 \wedge E \to E \wedge E \wedge E
$$
defines a diagonal map
$$
\Delta:E_\ast E \longr E_\ast E \otimes_{E_\ast} E_\ast E.
$$
With this map, the two units and the multiplication map 
$$
E_\ast E \otimes_{E_\ast} E_\ast E \longr E_\ast E
$$
the pair becomes a {\it Hopf algebroid}; that is, this pair represents a
a functor from graded rings to groupoids. One invariant of this
setup is the graded ring $R_E$ of {\it primitives} in $E_\ast$ 
defined by the equalizer diagram
$$
\xymatrix{
R_E \rto & E_\ast \ar@<.5ex>[r]^-{\eta_R} \ar@<-.5ex>[r]_-{\eta_L} &
E_\ast E.
}
$$
Note that if $R_E = E_\ast$, then $E_\ast E$ is a Hopf algebra
over $E_\ast$.

For any graded commutative ring, write $\Spec(R)$ for the covariant
functor on graded commutative rings it represents.
As a shorthand, I am going to say that
$$
\xymatrix{
\Spec(E_\ast E \otimes_{E_\ast} E_\ast E) 
\ar@<1.25ex>[r] \ar@<-1.25ex>[r]\ar[r]&
\ar@<.675ex>[l] \ar@<-.675ex>[l]\Spec(E_\ast E) \ar@<.75ex>[r] \ar@<-.75ex>[r] &
\ar[l]\Spec(E_\ast)
}
$$
is an affine groupoid scheme over $R_E$. (I leave you to label all
the maps and to check they satisfy the simplicial identities.)
This is an abuse of nomenclature, as we are
working in a grading setting, but it is nonetheless convenient:
many Hopf algebroids are best understood by describing the
functors they represent. If $E_\ast$ is $2$-periodic it is
possible to arrange matters so that this abuse goes away.
See Example \ref{2-period}.

A {\it comodule} over
this Hopf algebroid is a left $E_\ast$-module $M$ equipped with
a morphism of graded  $R_E$-modules
$$
\psi_M:M \longr E_\ast E \otimes_{E_\ast} M
$$
subject to the obvious compatibility maps with the structure
maps of the Hopf algebroid. If $X$ is a spectrum, then the
inclusion
$$
E \wedge \eta \wedge X:E \wedge S^0 \wedge X \longr E \wedge E \wedge X
$$
defines a comodule structure on $E_\ast X$. For example, if $X=S^0$,
then $E_\ast S^0 = E^\ast$ and
$$
\psi_{E_\ast} = \eta_R: E_\ast \longr E_\ast E.
$$
For much more
on comodules see \cite{Rav} \S A.1 and \cite{groupoid}.

Let's write $\comod_{E_\ast E}$ for the category of comodules over
the Hopf algebroid $(E_\ast,E_\ast E)$ and $\Mod_{E_\ast}$ for the modules
over $E_\ast$. 

\begin{rem}\label{ext-hopf-algebroid}The forgetful functor
$$
\comod_{E_\ast E} \to \Mod_{E_\ast}
$$
has a right adjoint $M \mapsto E_\ast E \otimes_{E_\ast} M$ -- the
{\it extended comodule} functor. Thus $\comod_{E_\ast E}$ has
enough injectives and we can define right derived functors; in particular,
$$
\Ext^s_{E_\ast E}(M,N)
$$
are the right derived functors of $\comod_{E_\ast E}(M,-)$. As
often happens, we can sometimes replace injectives by appropriate acyclic
objects. For example, if $M$
is {\it projective} as an $E_\ast$-module and
$$
N \to I^0 \to I^1 \to \cdots
$$
is a resolution by modules $I^s$ which are retracts of the extended
comodules, then a simple bicomplex argument shows that
$$
\Ext^s_{E_\ast E}(M,N) \cong H^s\comod_{E_\ast E}(M,I^\bullet).
$$
To be brutally specific, the triple resolution of $N$ obtained
from the composite functor $N \mapsto E_\ast E \otimes_{E_\ast} N$
yields a resolution with
$$
I^s = (E_\ast E)^{\otimes s+1} \otimes_{E_\ast} N.
$$
\end{rem}

\begin{rem}\label{alg-sphere}Of particular interest is the case of
$$
\Ext_{E_\ast E}^s(\Sigma^t E_\ast, N).
$$
Letting $t$ vary we have an equalizer diagram of $R_E$-modules
$$
\xymatrix{
\comod_{E_\ast E}(\Sigma^\ast E_\ast,N) \rto& 
N \ar@<.5ex>[r]^-{d_0} \ar@<-.5ex>[r]_-{d^1} &
E_\ast E \otimes_{E_\ast} N
}
$$
where $d^0(x) = 1\otimes x$ and $d^1(x) = \psi_N(x)$. This is equalizer
will be called the module of {\it primitives} and, hence, the functors
$\Ext^s_{E_\ast E}(\Sigma^\ast E_\ast,-)$ are the right derived functors
of the primitive element functor. The natural map
$$
\xymatrix@C=50pt{
M \rto^-{x \mapsto 1 \otimes x} & E_\ast E \otimes_{E_\ast E} M
}
$$
defines an isomorphism of $M$ to the primitives in the extended comodule;
hence, if we use the triple resolution of
the previous remark we have
$$
\Ext^s_{E_\ast E}(\Sigma^\ast E_\ast,N) = H^s(C^\bullet(N))
$$
where $C^\bullet(N)$ is the cobar complex. (See Remark \ref{bar-complex}
for formulas.) This normalization of this complex is the reduced
cobar complex $\bar{C}^\bullet(N)$.

If $A$ is a graded commutative algebra in $E_\ast E$-comodules,
then $C^\bullet (A)$ is a simplicial graded commutative $R_E$-algebra
and $\bar{C}^\bullet(A)$ is a differential graded commutative
algebra.
\end{rem}

We now define the Adams-Novikov Spectral Sequence (ANSS) for $E_\ast$.
Call a spectrum $K$ a {\it relative $E$-injective} if $K$ is a retract of a spectrum
of the form $E \wedge K_0$ for some $K_0$. The  following begins the
process. Compare Lemma \ref{basic}.

\begin{lem}\label{basic2}Let $X$ be a spectrum with $E_\ast X$
projective as an $E_\ast$-module. Then for all spectra $Y$
there is a Hurewicz map
$$
[X,Y] \longr \comod_{E_\ast E}(E_\ast X,E_\ast Y)
$$
which is an isomorphism if $Y$ is a relative $E$-injective.
\end{lem}

Then we have resolutions. Compare Definition \ref{adams-res}.

\begin{defn}\label{adams-res-e}Let $Y$ be a spectrum. Then
an $E_\ast$-Adams resolution for $Y$ is a sequence of spectra
$$
Y \to K^0 \to K^1 \to K^2 \to \cdots
$$
where each $K^s$ is a relative $E_\ast$-injective and for all 
relative $E_\ast$-injectives, the chain complex of
$R_E$-modules
$$
\cdots \to [K^2,K] \to [K^1,K] \to [K^0,K] \to [Y,K] \to 0
$$
is exact.
\end{defn}

Again such exist and are unique up to an appropriate notion of chain
equivalence. Furthermore, out of any $E_\ast$-Adams resolution we
get -- exactly as in Remark \ref{adams-tower} -- an Adams tower and
hence the Adams-Novikov Spectral Sequence
$$
\Ext_{E_\ast E}^s(\Sigma^t E_\ast X,E_\ast Y) \Longrightarrow
[\Sigma^{t-s}X,\lim Y_s].
$$
Again the spectrum $\lim Y_s$ is $E_\ast$-local and, under favorable
circumstances it is the $E_\ast$-localization $L_EY$ of $Y$. Convergence
remains an issue. Again see \cite{stlocal}. 
A sample nice result along these lines is the
following. The spectrum $E$ still satisfies the assumptions at the
beginning of this section.

\begin{lem}\label{first-conv}Let $E$ be a spectrum so that
\begin{enumerate}
\item $\pi_t E = 0$ for $t<0$ and either $\pi_0E \subseteq \QQ$ or
$\pi_0 E = \ZZ/n\ZZ$ for some $n$;

\item $E_\ast E$ is concentrated in even degrees.
\end{enumerate}
Suppose $Y$ is also has the property that $\pi_t Y = 0$ for
$t < 0$. Then 
\begin{enumerate}

\item the inverse limit of any $E_\ast$-Adams tower for $Y$ 
is the $E_\ast$-localization $L_E Y$;

\item $L_E Y = L_{H_\ast(-,\pi_0 E)}Y$;

\item $Ext^s(\Sigma^t  E^\ast,E_\ast Y) = 0$ for $s > t-s$ and the
ANSS converges.
\end{enumerate}
\end{lem}

\section{Complex oriented homology theories}

We begin to add to our assumptions on homology theories. In
the next definition we are going to use a natural isomorphism
$$
E^0 \cong \tilde{E}^2S^2
$$
where the target is the reduced cohomology of the two sphere.
There is an evident such isomorphism given by
$$
x \mapsto S^2 \wedge x:[S^0,E] \longr [S^2,S^2\wedge E],
$$
but we can always multiply by $-1$ to get another. We nail down the
isomorphism
by insisting that for $E = H_\ast(-,\ZZ)$ the isomorphism
should send $1\in H^0(\mathrm{pt},\ZZ)$ to the first Chern class of the canonical complex
line bundle over $S^2 = \CP^1$. We do this so that we
obtain the usual complex orientation for $H_\ast(-,\ZZ)$.

\begin{defn}\label{cmplx-oriented} Let $E_\ast$ be a multiplicative
homology theory. A {\it complex orientation} for $E^\ast$ is 
a {\bf natural} theory of Thom classes for complex vector bundles.
That is, given a complex $n$-plane bundle over a  $CW$-complex
$X$, there is a Thom class
$$
U_V\in E^{2n}(V,V-\{0\})
$$
so that
\begin{enumerate}
\item $U_V$ induces a Thom isomorphism
$$
\xymatrix{
E^k X \rto^\cong & E^kV \ar[rr]^-{U_V\smile} && E^{2n+k}(V,V-\{0\});
}
$$

\item if $f:X' \to X$ is continuous, then $U_{f^\ast V} = f^\ast U_V$.

\item if $\gamma$ is the canonical line over $\CP^1 = S^2$, then
the image of $U_\gamma$ under
$$
E^2(V(\gamma),V(\gamma)-\{0\}) \to E^2(V(\gamma)) \cong E^2(S^2)
$$
is the generator of $\tilde{E}^\ast (S^2)$ given by the image of
unit under the natural suspension isomorphism $E^0 \cong \tilde{E}^2S^2$.
\end{enumerate}
\end{defn}

If the homology theory $E_\ast$ is complex oriented, it  will have a good theory 
of Chern classes. In particular, any complex bundle $V$ has an
Euler class $e(V)$ (the top Chern class) given as the image of $U_V$ under
the map
$$
E^{2n}(V,V-\{0\}) \to E^{2n}(V) \cong E^{2n}X.
$$
From this and the Thom isomorphism one gets a Gysin sequence and
the usual argument now implies that 
\begin{equation}\label{ecpinfty}
E^\ast \CP^\infty \cong E^\ast[[e(\gamma)]]
\end{equation}
where $\gamma$ is the universal $1$-plane bundle. In fact,
$$
E^\ast((\CP^\infty)^n) \cong E^\ast[[e(\gamma_1),\ldots,e(\gamma_n)]]
$$
where $\gamma_i = p_i^\ast \gamma$ is the pull-back of $\gamma$
under the projection onto the $i$th factor. The method of universal
examples now implies that there is a power series
$$
x+_Fy \defeq F(x,y) \in E^\ast[[x,y]]
$$
for calculating the Euler class of the tensor product of two line bundles:
$$
e(\xi \otimes \zeta) = e(\xi) +_F e(\zeta).
$$
The usual properties of the tensor product imply that the power series
$x+_Fy$ satisfies the following properties
\begin{enumerate}

\item (Unit) $x+_F0 = x$;

\item (Commutativity) $x+_Fy = y+_Fx$; and

\item (Associativity) $(x+_Fy)+_Fz = x+_F(y+_Fz)$.
\end{enumerate}
A power series satisfying these axioms is called a {\it formal group law};
thus a choice of complex orientation for $E_\ast$ yields a formal group
law over $E^\ast$.\footnote{This is a smooth one-dimensional 
formal group law, but this is the only kind we'll have, so the adjectives
aren't necessary. Also, if you're into algebraic geometry, you realize that
the important object here is not $F$, but the formal group scheme represented
by $E^\ast \CP^\infty$ -- which does not depend on the choice of orientation.
Such formal group schemes have been called by one-parameter
formal Lie groups. See \cite {smith} for the background on nomenclature.
Here we're doing computations here, so we'll stick to the formulas.}

Now suppose we have two different orientations for $E_\ast$ and hence
two Euler classes $e(V)$ and $e'(V)$. An examination of Equation
\ref{ecpinfty} shows that there must be a power series $f(x) \in E^\ast[[x]]$
so that
$$
e'(\gamma) = f(e(\gamma)).
$$
The normalization condition of Definition \ref{cmplx-oriented}.3 implies
that $f'(0)=1$; thus, $f$ is an element of the group of power series
invertible under composition. Again checking universal examples
we see that if $F$ and $F'$ are the two resulting formal group laws,
then we see that
$$
f(x+_Fy) = f(x) +_{F'} f(y).
$$
Such a power series is called a {\it strict isomorphism} of formal group
laws.

\begin{exam}[{\bf Complex cobordism}] Formal group laws and their
isomorphisms assemble into a functor $\fgl$ to groupoids from the category
of graded rings: given a graded ring $R$ the objects of $\fgl(R)$ are
the formal group laws over $R$ and the morphisms are the strict
isomorphisms of formal group laws. If $g:R \to S$ is a homomorphism
of graded commutative rings, we will write
$$
f^\ast:\fgl(R) \longr \fgl(S)
$$
for the resulting morphism of groupoids.

It is easy to see that this is
an affine groupoid scheme (in the graded sense) over $\ZZ$; a theorem
of Lazard calculates the rings involved. The key observation  -- which
makes smooth, 1-dimensional formal group laws special -- is the symmetric
$2$-cocyle lemma. This says that given a formal group law $F$ over
a graded ring $R$, there are elements $a_1,a_2,\cdots$ of $R$
and an equality
$$
x+_Fy = x+y + a_1C_2(x,y) + a_2C_3(x,y)+\cdots
$$
in $R/(a_1,a_2,\ldots)^2$ where $C_n(x,y)$ is the $n$th homogeneous
symmetric $2$-cocyle
$$
C_n(x,y) = (1/d)[(x+y)^n - x^n - y^n].
$$
Here $d=p$ if $n$ is a power of a prime $p$; otherwise $d=1$. It follows
easily from this that objects of $\fgl(-)$ are represented by the Lazard
ring, which is  a polynomial ring
$$
L \cong \ZZ[a_1,a_2,\ldots]
$$
Note this isomorphism is canonical only modulo $(a_1,a_2,\ldots)^2$. The morphisms in $\fgl(-)$ are represented by
$$
W=L[t_1,t_2,\ldots]
$$
where the $t_i$ are the coefficients of the universal strict isomorphism.
The elements $a_i$ and $t_i$ have degree $2i$ (in the homology grading).
The pair $(L,W)$ is a Hopf algebroid over $\ZZ$.

A celebrated theorem of Quillen \cite{quillen}
implies that any choice of a complex
orientation for the complex cobordism spectrum $MU$ yields
an isomorphism of Hopf algebroids
$$
(MU_\ast,MU_\ast MU) \cong (L,W).
$$
\end{exam}

\def\ZZP{{{\mathbb{Z}_{{(p)}}}}}

\begin{exam}[{\bf Brown-Peterson theory}]\label{bptheory}The integers
$\ZZ$ has some drawbacks as a ring -- most dauntingly, it is not a local ring.
We often prefer to work over $\ZZP$-algebras for some fixed prime
$p$. Then we can cut down the formal group laws in question.

For any integer $n$ and any formal group law over any ring $R$ we
can form the $n$-series
$$
[n]_F(x) = [n](x) = x+_F + \cdots +_F x = nx + \cdots 
$$
where the sum is taken $n$-times. If $R$ is a $\ZZP$-algebra and 
$n$ is prime to $p$, this power series is invertible (under composition);
let $[1/n](x)$ be the inverse.  A formal group law is $p$-typical if for all
 $n$ prime to $p$ the formal sum
$$
[1/n](x +_F \zeta x +_F \zeta^2 x +_F \cdots +_F \zeta^{n-1} x) = 0.
$$
Here $\zeta$ is a primitive $n$th root of unity. The sum is {\it a priori}
over $R[\zeta]$; however, since it is invariant under the Galois action, it is actually
over $R$.

This definition, admittedly rather technical, has the following implications
(see \cite{Rav} \S A.2 for pointers to original sources):
\begin{enumerate}

\item Let $F$ be any formal group law over a $\ZZP$-algebra $R$. Then
there is a strict isomorphism
$$
e_F:F \to G
$$
from $F$ to a $p$-typical formal group law $G$. If $F$ is already
$p$-typical, then $e_F$ is the identify. 

\item The idempotent $e$ is natural in $F$; that is, if $\phi:F \to F'$ is a strict isomorphism
of formal group laws, then there is a unique isomorphism $e_{\phi}$ making
the following diagram commute\footnote{This statement is obvious,
since all the maps are isomorphisms. I added this to belabor the
functoriality of the idempotent.}
$$
\xymatrix{
F \rto^{e_F} \dto_\phi & G \dto^{e_\phi}\\
F' \rto_{e_{F'}} & G'.
}
$$

\item For any $p$-typical formal group law $F$ over a $\ZZP$-algebra $R$, 
there are elements $v_i \in R$ so that
\begin{equation}\label{p-typical-p-series}
[p]_F(x) = px +_F v_1x^p +_F v_2x^{p^2} +_F \cdots.
\end{equation}

\item Notice that we could write $v_n(F)$ for the elements $v_n \in R$
of Equation \ref{p-typical-p-series} as they depend on the formal
group law $F$. Now we can prove that if
$F$ and $G$ are two $p$-typical formal group
laws over $R$ so that $v_n(F) = v_n(G)$ for all $n$, then $F=G$. We can also prove
that given element $a_n$ of degree $2(p^n-1)$, there is a
$p$-typical formal group $F$ with $v_n(F) = a_n$.

\item If $\phi:F \to F'$ is any isomorphism of $p$-typical formal group
laws over a $\ZZP$-algebra $R$, then there are elements $t_i \in R$
so that
$$
\phi^{-1}(x) = x +_F t_1x^p +_F t_2x^{p^2} +_F \cdots.
$$
\end{enumerate}

Let $\fglp$ denote the functor from graded $\ZZP$-algebras to
groupoids which assigns to any $R$ the groupoid of $p$-typical
formal group laws over $R$ and all their isomorphisms. Then
points (1) and (2) show that there is a natural equivalence of
groupoids (the {\it Cartier idempotent})
$$
e:\fgl \to \fglp
$$
which is idempotent. Points (3), (4), and (5) show that $\fglp$ is
represented by the Hopf algebroid $(V,VT)$ where
$$
V = \ZZP[v_1,v_2,\ldots]
$$
where the $v_i$ arise from the universal $p$-typical formal group
as in point (3)  and
$$
VT = V[t_1,t_2,\ldots]
$$
where the $t_i$ arise from the universal isomorphism as in point (5).
The degrees of $v_i$ and $t_i$ are $2(p^i-1)$. Note that these degrees
grow exponentially in $i$, which often helps in computations.

Quillen \cite{quillen} noticed that there is a complex oriented ring spectrum $BP$
(constructed, without its ring structure, by Brown and Peterson) and an
isomorphism of Hopf algebroids
$$
(BP_\ast,BP_\ast BP) = (V,VT).
$$
\end{exam}

\begin{exam}[{\bf Landweber exact theories}]Let $F$ be a formal
group law over a graded ring $A$ classified by map of graded
rings $F:L \to A$. We would like the functor
$$
X \mapsto A \otimes_L MU_\ast X
$$
to define a homology theory $E(A,F)$. Using Brown representability,
this will happen if the functor $A \otimes_L (-)$ is
exact on the category of $(MU_\ast,MU_\ast MU)\cong (L,W)$-comodules.
The Landweber Exact Functor Theorem (LEFT) (see \cite{landleft} and
\cite{LEFT}) writes down a condition
for this to be true.

We can write down a formal requirement fairly easily. Assuming for the
moment that $E(A,F)$ exists, we have that
$$
E(A,F)_\ast E(A,F) \cong A \otimes_L MU_\ast E(A,F)
$$
and $MU_\ast E(A,F) = MU_\ast MU \otimes_L A$, by switching the
factors. Thus
\begin{equation}\label{land-algebroid}
E(A,F)_\ast E(A,F) \cong A \otimes_L W \otimes_L A.
\end{equation}
Gerd Laures has noticed that $A \otimes_L (-)$ is exact on
$(L,W)$ comodules if the {\it right} inclusion
$$
\eta_R:L \longr A \otimes_L W
$$
remains flat. For a computational condition, it's easier to work at a prime
$p$ and suppose that $F$ is $p$-typical. Then Landweber's criterion
is that the sequence $p,v_1,v_2,\cdots$ is a regular sequence
in $A$. The $p$-typical requirement can be removed:
it turns out that $v_i$ is defined for all formal group laws modulo
$p,v_1,\ldots,v_{i-1}$ and we simply require that the resulting
sequences be regular for all $p$. If $A$ is $\ZZP$-algebra and $q$
is prime to $p$ this is automatic for the prime $q$.

If the formal group $F$ satisfies these conditions, then we say
$E(A,F)$ is a Landweber exact theory. The equation \ref{land-algebroid}
now implies that there is an isomorphism of Hopf algebroids
$$
(E(A,F)_\ast, E(A,F)_\ast E(A,F)) \cong (A, A \otimes_L W \otimes_L A).
$$
This represents the (graded) affine groupoid scheme which assigns to
any ring $R$ the groupoid with objects the ring homomorphisms
$f:A \to R$. An isomorphism $\phi:f \to g$ will be a strict isomorphism
$$
\phi:f^\ast F \to g^\ast F
$$
of formal group laws.

Important first examples of Landweber exact theories are the
{\it Johnson-Wilson} theories $E(n)$ with
$$
E(n)_\ast \cong \ZZP[v_1,\cdots,v_{n-1},v_n^{\pm 1}].
$$
We normally choose a $p$-typical formal group law $F$ over $E(n)_\ast$
with $p$-series
$$
[p](x) = px +_F v_1x^p +_F \cdots +_F v_nx^{p^n}.
$$
Since $E(n)$ is not connected, the localization $L_{E(n)} = L_n$ can
be a radical process. Nonetheless we have the Hopkins-Ravenel
{\it chromatic convergence} theorem \cite{orange}
which says that there are natural
transformation of functors $L_n \to L_{n-1}$ and, if $X$ is a finite
$CW$ spectrum, then the induced map
$$
X \longr \holim L_n X
$$
is localization at $H_\ast(-,\ZZP)$. See \cite{orange}.
\end{exam}

\begin{exam}[{\bf K-theory}]\label{k-theory} Let $K$ be complex
$K$-theory. The Euler class of a complex line bundle is normally
defined to be
$$
e'(\gamma) = \gamma - 1 \in \tilde{K}^0(X).
$$
With this definition, 
$$
e'(\gamma_1 \otimes \gamma_2) = e'(\gamma_1) + e'(\gamma_2)
+ e'(\gamma_1)e'(\gamma_2).
$$
The resulting formal group law is the multiplicative formal group law
$$
x+_{\GG_m} y = x + y + xy \in K^0[[x,y]].
$$
Since we have arranged our definitions so that the Euler class lies
in $K^2(X)$, we redefine the Euler class so that
$$
e(\gamma) = u^{-1}e'(\gamma) = u^{-1}(\gamma - 1)
$$
where $u \in K_2 = K^{-2}$ is the Bott element. The graded formal
group law then becomes
$$
x +_{F_u} y = x+y+uxy = u^{-1}(u x + {\GG_m} uy).
$$
This has $p$-series $(1+ux)^p -1$; hence $K$-theory is
Landweber exact. (But $\ZZP \otimes K_\ast(-)$ is not $p$-typical.)
The groupoids scheme represented by the pair $(K_\ast,K_\ast K)$
assigns to each graded ring $R$ the groupoid with objects
 the units $v \in R_2$ and morphisms the strict isomorphisms
 $$
 \phi:F_v \to F_{w}.
 $$
 By replacing $\phi(x)$ by $w^{-1}\phi(vx)$ we see that this is
 isomorphic to the groupoid with objects the units $v$ and
 morphisms the non-strict isomorphisms $\phi$ of $\GG_m$ so
 that $\phi'(0) = w^{-1}v$.
 \end{exam}
  
\begin{exam}[{\bf 2-periodic theories}]\label{2-period} With the
example of $K$-theory (see \ref{k-theory})
we first encounter a $2$-periodic theory. Other important
examples are the Lubin-Tate theories $E_n$.  See \ref{grp-coh-2}
below.

Let $A$ be
a graded ring of the form $A = A_0[u^{\pm 1}]$ where $A_0$ is
in degree $0$ and $u$ is a unit in degree $2$. Let $F$ be a formal
group law over $A$. Then we can define an ungraded formal group
law $G$ over $A_0$ by the formula
$$
x +_G y = u( u^{-1}x +_F u^{-1} y).
$$
Let $\Ga = A \otimes_L W \otimes_L A$.
The groupoid scheme represented by $(A,\Ga)$
assigns to each graded commutative ring the groupoid with objects
$f:A \to R$ and morphism the strict isomorphisms $\phi:f_0^\ast F \to
f_1^\ast F$. This is isomorphic to the groupoid with objects $(g,v)$
where $g:A_0 \to R_0$ is a ring homomorphism and $v \in R_2$ is
a unit. The morphisms are all isomorphisms $\phi:g_0^\ast G
\to g_1^\ast G$ so that $v_1 = \phi'(0)v_0$. Thus if $a_i \in \Ga_{2i}$
are the coefficients of the universal strict isomorphism from
$\eta_L^\ast F$ to $\eta_R^\ast F$, then
$$
b_i = \brackets{u^{-i}a_i,}{i> 0;}{u^{-1}\eta_R(u),}{i=0}
$$
are the coefficients of the universal isomorphism from 
$\eta_L^\ast G$ to $\eta_R^\ast G$.

There is an equivalence of categories between $(A,\Ga)$-comodules
concentrated in even degrees\footnote{There is
a theory that includes odd degrees; this
is assumption is for convenience. See \cite{LEFT}} and $(A_0,\Ga_0)$-comodules,
where $(-)_0$ means the elements in degree $0$. The functor
one way simply send $M$ to $M_0$. To get the functor back,
let $N$ be an $(A_0,\Ga_0)$-comodule and, define, for $i$, new
$(A_0,\Ga_0)$-comodule $M(i)$ with $M(i) = u^iM$ and
$$
\psi_{M(i)}(u^ix) = b_0^i u^i\psi_M(x).
$$
If $M$ is an $(A,\Ga)$-comodule in even degrees, then $M = \oplus M_0(i)$,
where $M_0(i)$ is now in degree $2i$. 

Using this equivalence of categories we have that
$$
\Ext_\Ga(A,M) \cong \Ext_{\Ga_0}(A_0,M_0)
$$
and if $M$ is concentrated in even degrees,
$$
\Ext_\Ga(\Sigma^{2t}A,M) \cong \Ext_\Ga(A,\Sigma^{-2t}M)
\cong \Ext_{\Ga_0}(A_0,M_0(t)).
$$
Also $\Ext_\Ga(\Sigma^{2t+1}A,M)=0$.
\end{exam}

\begin{exam}[{\bf The Morava stabilizer group}]\label{Morava-stabilizer}Let $F$ be a formal
 group law over a graded commutative ring $R$; it is possible
 for the Hopf algebroid to be important even if the pair $(R,F)$
 is not Landweber exact. Here is a crucial example. Set
 $$
 K(n)_\ast = \FF_p[v_n^{\pm 1}]
 $$
 and let $F_n$ be the unique $p$-typical formal group law with $p$-series
 $[p](x) = v_nx^{p^n}$. This is the {\it Honda} formal group law. Then
 the Hopf algebroid 
 $$
 (K(n)_\ast,K(n)_\ast \otimes_L W \otimes_L K(n)_\ast)
 $$
 If we write 
 $$
 \Sigma(n) \defeq K(n)_\ast \otimes_L W \otimes_L K(n)_\ast
 $$
 then $\Sigma(n)$ is a graded Hopf {\it algebra} over
 $K(n)_\ast$ -- since $v_n$ is primitive -- by looking at the
 $p$-series we have
 $$
 \Sigma(n) \cong \FF_p[v_n^{-1}][t_1,t_2,\cdots]/(t_i^{p^n}-v_n^{p^i-1}t_i)
 $$
 is a graded Hopf algebra.
  
 This Hopf algebra can be further understood as follows. Define a faithfully 
 flat extension
 $$
 i:\FF_p[v_n^{\pm 1}] \longr \FF_{p}[u^{\pm 1}]
 $$
 where $u^{p^n-1} = v_n$. Because $\FF_{p^n}[u^{\pm 1}]$
 is now $2$-periodic, we can (as in the last example) rewrite
 the Honda formal group law $i^\ast F_n$ as an (ungraded)
 formal group law $\Ga_n$ over
 $\FF_{p}$ with $p$-series exactly $x^{p^n}$. Then
 $$
 \FF_p[u^{\pm 1}] \otimes_{K(n)_\ast} (K(n)_\ast,\Sigma(n))
 \cong \FF_p[u^{\pm 1}] \otimes (\FF_p, S(n))
 $$
 where $S(n)$ is a  Hopf algebra concentrated in degree $0$
 which can be written
 $$
 S(n) = \FF_p[x_0^{\pm 1},x_1,\cdots]/(x_i^{p^n}-x_i).
 $$
 Then $\Spec(S(n))$
 is the group scheme which assigns to each commutative $\FF_p$-algebra
 $R$ the (non-strict) automorphisms of $\Ga_n$ over $R$.
 The $\bar{\FF}_p$ points  of this group scheme -- 
 the automorphisms of $\Ga_n$ over the algebraic closure
 of $\FF_p$ -- is called the {\it Morava stabilizer group} $S_n$. From the
 form of $S(n)$ or from the fact that the $p$-series of $\Ga_n$ is simply
 $x^{p^n}$, we see that all these automorphisms are realized over
 $\FF_{p^n}$. We take this group up below in Section \ref{sect-group}.
 For now I note that it is known and has been heavily studied: 
 see \cite{Rav}\S A2.2 and \cite{hwhduke} -- and from either of these
 references you can find many primary sources.
 If $n=1$ it is the group $\ZZ_p^\times$
 of  units in the $p$-adic integers; for $n \geq 2$ it is non-abelian. 
 For a thorough review of these groups, their cohomology and
 many more references, see \cite{hwhduke}.
 \end{exam}
 
\section{The height filtration}

Let $f:F \to G$ be a homomorphism of formal group laws over
a commutative ring $R$; thus
$$
f(F(x,y)) = G(f(x),f(y)).
$$
Differentiating with respect to $y$ and setting $y=0$, we get 
a formula
\begin{equation}\label{inv-diif-form}
f'(x)F_y(x,0) = G_y(f(x),0)f'(0).
\end{equation}
Here $F_y(x,y)$ is the partial derivative of $F(x,y)$ with respect to $y$.
The power series $F_y(x,0) = 1 + \hbox{higher terms}$ and is therefore
invertible. Thus if $f'(0) = 0$, then $f'(x) = 0$. If $R$ is a $\QQ$-algebra,
this implies $f(x) = 0$; however, if $R$ is an $\FF_p$ algebra,
this simply implies that there is a power series $g(x)$ so that
$$
f(x) = g(x^p).
$$
If $R$ is an $\FF_p$ algebra, we let $\sigma:R \to R$ denote the 
Frobenius, and if $F$ is a formal group law over $R$, we let
$$
F^{(p)} = \sigma^\ast F
$$
be the pull-back of $F$ along the Frobenius. Thus if
$$
F(x,y) = \sum a_{ij}x^iy^j,
$$
then 
$$
F^{(p)}(x,y) = \sum a_{ij}^px^iy^j.
$$
The power series $\sigma(x) = x^p$ defines a homomorphism
$$
\sigma:F \to F^{(p)}.
$$
Combining all these remarks, we have the following:

\begin{lem}\label{frob-factor}Let $f:F \to G$ be a homomorphism of
formal group laws over an $\FF_p$-algebra $R$ and suppose
that $f'(0) = 0$. Then there is a unique homomorphism
of formal group laws $g:F^{(p)} \to G$ and a factoring
$$
\xymatrix{
F \rto^{\sigma} \ar@/_1pc/[rr]_f& F^{(p)} \rto^g &G.
}
$$
\end{lem}

This can be repeated if $g'(0) = 0$ to get a factoring
$$
\xymatrix{
F \rto^\sigma \ar[drr]_f& F^{(p)} \rto^{\sigma} & F^{(p^2)} \dto^{g_2}\\
&& G
}
$$

\begin{defn}\label{height}Let $F$ be a formal group law over 
an $\FF_p$-algebra $R$. Then $F$ has a {\bf height} at least $n$
if there is a factoring of the $p$-series
$$
\xymatrix{
F \rto^\sigma \ar[drrr]_{[p]}& F^{(p)} \rto^{\sigma} & \cdots \rto^\sigma &F^{(p^n)} \dto^{g_n}\\
&&& F
}
$$
\end{defn}

Every formal group law over an $\FF_p$-algebra has height at least
$1$. The multiplicative formal group law $\GG_m (x,y) = x+y+xy$
has height $1$ but not height $2$; the additive formal group law
$\GG_a(x,y) = x+y$ has height at least $n$ for all $n$; thus,
we say it has infinite height. The Honda formal group law
of \ref{Morava-stabilizer} has height exactly $n$.

By construction, if $F$ has height at least $n$ there is a power
series $g_n(x) = v^F_nx + \cdots$ so that
\begin{equation}\label{p-series-n}
[p]_F(x) = g_n(x^{p^n}) = v^F_nx^{p^n} + \cdots
\end{equation}
and $F$ has height at least $n+1$ if and only if $v^F_n=0$. The element
$v^F_n$ depends only on the strict isomorphism class of $F$.
Thus, to check whether $v_n = 0$, we can assume that $F$
is $p$-typical, where we already know that
$$
[p]_F(x) = px +_F v_1 x^p +_F \cdots +_F v_n x^{p^n} + \cdots
$$
This formula implies that the two definitions of $v_n$-agree. It
also implies that there is an ideal
$$
I_n \defeq (p,\ldots,v_{n-1}) \subseteq BP_\ast
$$
so that a $p$-typical formal group law $F$ has height at least $n$
if and only if $0=FI_n \subseteq R$, where $F:BP_\ast \to R$ classifies
$F$. The image of $I_n$ in $R$ depends only on the strict isomorphism
class of $F$; from this is follows that the ideal $I_n$ is {\it invariant} -- that is, 
it is a subcomodule of $BP_\ast$. It is a theorem of Landweber's (see
\cite{invar}) that
the $I_n$ and $\{0\}$ are  the only prime (and, in fact, radical) invariant
ideals of $BP_\ast$.

We now have a filtration of the groupoid scheme $\fglp$ over
$\ZZP$ the groupoid schemes
$$
\cdots \subseteq \fglp(n) \subseteq \fglp(n-1)\subseteq \cdots \subseteq \fglp(1) \subseteq 
\fglp
$$
where $\fglp(n)$ assigns to each commutative ring $R$ the groupoid
of $p$-typical formal groups $F$ for which $FI_n = 0$. Since these
are defined by the vanishing of an ideal, we call these closed; the
open complement of $\fglp(n)$ is the subgroupoid scheme
$\cU(n-1) \subseteq \fglp$
of $p$-typical formal groups $F$ for which $FI_n = R$. (Note the shift
in numbering.) This gives a filtration
$$
\cU(0) \subseteq \cU(1) \subseteq \cdots \cU(n-1) \subseteq \cU(n) \subseteq
\cdots \subseteq 
\fglp.
$$
Note that neither filtration is exhaustive as
$$
\cap \fglp(n) = \fglp - [\cup\ \cU(n-1)]
$$
contains the additive formal group law over $\FF_p$.

\begin{rem}\label{filt-basic}
Here are some basic facts about this filtration:
\begin{enumerate}

\item The groupoid scheme $\fglp(n)$ is defined by the Hopf algebroid
\begin{align*}
(BP_\ast/I_n, BP_\ast/I_n &\otimes_{BP_\ast} BP_\ast BP)\\
&\cong (\FF[v_n,v_{n+1},\ldots],\FF[v_n,v_{n+1},\ldots][t_1,t_2,\ldots]).
\end{align*}

\item A $p$-typical formal group law $F$ over $R$ lies in $\cU(0)$
if and only if $R$ is a $\QQ$-algebra; $F$ is then isomorphic
to the additive formal group $\GG_a$ over $R$ (induced from $\QQ$).
The automorphisms of $\GG_a$ over a $\QQ$-algebra $R$
is isomorphic to the group of units $R^\times$.

\item The groupoid scheme $\cU(n)$, $n \geq 1$, is not isomorphic
or even equivalent to an affine groupoid scheme given by
a Hopf algebroid. It is, however, equivalent
{\it locally for the flat topology} (see Definition \ref{local-fl-top} for this notion.)
to the groupoid scheme of defined
by the Hopf algebroid of the Johnson-Wilson theories
$(E(n)_\ast,E(n)_\ast E(n))$. This is a result of Hovey-Strickland
\cite{HSt} and Naumann \cite{naumann}.

\item The relative open $\fglp(n) \cap \cU(n)$ consists of those
$p$-typical formal group laws $F$ over an $\FF_p$ algebra $R$
which are of {\it exact} height $n$; that is, $FI_n = 0$
and $v^F_n \in R$ is a unit. This groupoid scheme has 
exactly one ``geometric point''; that is, over a algebraically
closed field there is a unique isomorphism class of formal
groups of height $n$. Its automorphism group was discussed below
in Section \ref{sect-group}.
\end{enumerate}
All of this should be treated rigorously with the language of stacks.
See \cite{smith}.
\end{rem}

\def\cX{{{\mathcal{X}}}}

\begin{rem}[{\bf Quasi-coherent sheaves over groupoid schemes}]\label{qc}
Because groupoid scheme $\cU(n)$ is not isomorphic to the groupoid
scheme given by a Hopf algebroid; thus the scope of inquiry must
expand a little at this point. In particular, we need to find a reformulation
of the notion of comodule which can be adapted to such groupoid
scheme. These are the {\it quasi-coherent sheaves}. For much
more on this topic see \cite{groupoid} or, for  the general
theory \cite{LM}

Let $\cX$ be a groupoid scheme. If $\cX$ is not affine, we have
to take some care about what this means. At the very least it
includes the provision that
the objects and morphisms form sheaves in the flat topology.
Thus, if $\cX_0$ is the objects functor of $\cX$,
then for all faithfully flat extensions 
$R \to S$ there is an equalizer diagram 
$$
\xymatrix{
\cX_0(R) \rto & \cX_0(S) \ar@<.5ex>[r] \ar@<-.5ex>[r] &
\cX_0(S \otimes_R S)
}
$$
There is a similar condition on morphisms. The examples $\cU(n)$
above satisfy this conditions, as does an groupoid scheme arising
from a Hopf algebroid.

If $R$ is a ring, write $\Spec(R)$ for the functor it represents. (Remember,
we are working with graded rings, for the most part, so this is a bit of
an abuse of notation.) The functor $\Spec(R)$ is trivially a groupoid
scheme, with only identity morphisms. Then the morphisms
of groupoid schemes
$$
x:\Spec(R) \to \cX
$$
are in one-to-one correspondence with objects $x\in \cX(R)$.
A $2$-commuting diagram 
\begin{equation}\label{2-commute}
\xymatrix@R=10pt{
\Spec(S) \ar[dd]_{\Spec(f)} \ar[dr]^y\\
&\cX\\
\Spec(R) \ar[ur]_x
}
\end{equation}
is a pair $(f,\phi)$ where $f:R \to S$ is a ring homomorphism 
and $\phi:y \to f^\ast x$ is an isomorphism in $\cX(S)$.

A quasi-coherent sheaf  $\cF$ assigns to each morphism of
groupoid schemes $x:\Spec(R) \to \cX$ an $R$-module
$\cF(R,x)$ and to each $2$-commuting diagram \ref{2-commute}
a morphism of $R$-modules $\cF(R,x) \to \cF(S,y)$ so that the
induced map 
$$
S\otimes_R \cF(R,x) \longr \cF(S,y)
$$
is an isomorphism.

If $\cX$ is the groupoid scheme arising from the Hopf algebroid
$(A,\Ga)$, then an object $x\in \cX(R)$ is a homomorphism
$x:A \to R$. If $M$ is a comodule we set $\cF(R,x) = R\otimes_A M$.
This is part of an equivalence of categories.

If two groupoid schemes are locally equivalent in the flat
topology, they have equivalent categories of quasi-coherent
sheaves; this implies they have equivalent sheaf cohomologies
for these sheaves (i.e., coherent cohomology). If $\cX$ arises
from a Hopf algebroid $(A,\Ga)$ this cohomology is exactly
$\Ext_\Ga^\ast(\Sigma^\ast A,-)$. See Theorem \ref{change-of-rings}.
\end{rem}

\begin{rem}[{\bf The invariant differential}]\label{inv-diff}Let $F(x,y)$
be a formal group law over $R$. The power series $F_y(x,0)$ introduced
at the beginning of this section didn't come from the ether. We define
a differential $f(x) \in R[[x]]dx$ to be invariant if
$$
f(x+_F y)d(x+_Fy) = f(x)dx + f(y)dy.
$$
Expanding this formula out and setting $y=0$, we see that the
invariant differentials form a free
sub-$R$-module of rank $1$ of $R[[x]]dx$  generated by
$$
\omega_F = \frac{dx}{F_y(x,0)};
$$
indeed, if $f(x)dx$ is an invariant differential, then
$f(x)dx = f(0)\omega_F.$
If $f:F\to G$ is a homomorphism of formal groups then the formula
of \ref{inv-diif-form} says that 
$$
f^\ast \omega_G = f'(0)\omega_F.
$$
If $R$ is a $\QQ$-algebra, then the integral of $\omega_F$ with zero
constant term defines
an isomorphism of $F$ with the additive formal group law.
\end{rem}

\section{The chromatic decomposition}

The chromatic decomposition attempts to isolate what part
of the Ext groups $\Ext_{BP_\ast BP}^\ast (\Sigma^\ast BP_\ast, BP_\ast)$
arise from the formal groups of exact height $n$. The result
is a spectral sequence.

To begin, we can easily isolate the part of the Ext groups arising
from formal groups over the rational numbers (exact height $0$
if you want) by inverting $p$. Thus we have a short exact
sequence of comodules
$$
0 \to BP\ast \to p^{-1}BP_\ast \to BP_\ast/(p^\infty) \to 0.
$$
The comodule $BP_\ast/(p^\infty)$, because it $p$-torsion,
lives (in a sense I won't make precise here) on $\fglp(1)$.
We'd like something that lives on $\fglp(1) \cap \cU(1)$ -- on
groupoid scheme of formal groups of exact height $1$. To
do this we'd like to form $v_1^{-1}BP_\ast/(p^\infty)$. For this
we need the following result.

\begin{lem}\label{invert-vn}Let $M$ be $(BP_\ast,BP_\ast BP)$-comodule
which is $I_n$-torsion as a $BP_\ast$-module. Then there is
a unique comodule structure on $v_n^{-1}M$ so that the
natural homomorphism $M \to v_n^{-1}M$ is a comodule map.
\end{lem}

\begin{proof}Write $\eta_R(v_n) = v_n - f$ where 
$$
 f \in I_n \otimes_{BP_\ast} BP_\ast BP\cong BP_\ast BP_\ast \otimes_{BP_\ast}I_n.
$$
Suppose we have
a comodule structure map $\psi$ for $v_n^{-1}M$ extending
$\psi_M$. Then for
all $a \in v_n^{-1}M$, there is a $k$ so that $v_n^ka \in M$.
Thus we have
$$
\psi_M(v_n^ka) = \eta_R(v_n^k)\psi(a)
$$
or
\begin{align*}
\psi(a) &= \eta_R(v_n)^{-k}\psi_M(v_n^ka)\\
&= v_n^{-k}(1-f/v_n)^{-k}\psi_M(v_n^ka).
\end{align*}
Expanding $(1-f/v_n)^{-k}$ as a power series and using that
$M$ is $I_n$-torsion we see that this equation defines $\psi$ and
proves its uniqueness.
\end{proof}
We now form the following diagram, where each of the triangles
$$
\xymatrix{
K \ar[d] & \ar@{-->}[l] N\\
M \ar[ur]
}
$$
represents a short exact sequence 
$$
0 \to K \to M \to N \to 0
$$
in comodules. The dotted arrow indicates that we will have
a long exact sequence in $\Ext$-groups. Using Lemma
\ref{invert-vn} we now inductively construct a diagram 
extending infinitely far to the right:
$$
\xymatrix@C=10pt{
BP_\ast \ar[d] & \ar@{-->}[l] BP_\ast/(p^\infty) \ar[d] &
\ar@{-->}[l] BP_\ast/(p^\infty,v_1^\infty)\ar[d]&
\ar@{-->}[l] BP_\ast/(p^\infty,v_1^\infty,v_2^\infty)\\
p^{-1}BP_\ast \ar[ur]&v_1^{-1}BP_\ast/(p^\infty) \ar[ur]
&v_2^{-1}BP_\ast/(p^\infty,v_1^\infty) \ar[ur]
}
$$
Let's write $BP_\ast/I_n^\infty$ for $BP_\ast/(p^\infty,\ldots,v_{n-1}^\infty)$.
We now apply the functor
$$
H^{s,t}(-) \cong \Ext^s_{BP_\ast BP}(\Sigma^t BP_\ast,-)
$$
to these short exact sequences to get an exact couple; the dotted
arrows now become maps $H^{s,t}(N) \to H^{s+1,t}(K)$. The
resulting spectral sequence -- {\it the chromatic
spectral sequence} --  then reads
\begin{equation}\label{chrom-ss}
E_1^{p,q}=H^{q,\ast}(v_p^{-1}BP/I_p^\infty) \Longrightarrow 
H^{p+q}(BP_\ast)
\end{equation}
with differentials
$$
d_r:E_r^{p,q} \to E_r^{p+r,q-r+1}.
$$
In particular, $E_2$ is the cohomology of the very interesting chain
complex
$$
0 \to H^{\ast,\ast}(p^{-1}BP_\ast) \to H^{\ast,\ast}(v_1^{-1}BP_\ast/I_1^\infty)
\to H^{\ast,\ast}(v_1^{-2}BP_\ast/I_2^\infty) \to \cdots
$$
and there is an edge homomorphism
$$
\xymatrix{
E_2^{p,0} \ar@{>>}[r] & E_\infty^{p,0}\ \ \ar@{>>->}[r]& \Ext_{BP_\ast BP}^p(\Sigma^\ast BP_\ast,BP_\ast)
}
$$
with $E_2^{p,0}$ a subquotient of
$$
Ext_{BP_\ast BP}^p(\Sigma^\ast BP_\ast,v_p^{-1}BP_\ast/I_p^\infty).
$$

\begin{exam}[{\bf The Greek letter elements}]It is easy to produce
an interesting family of permanent cycles in this spectral
sequence. There is a canonical inclusion of comodules
$$
\iota_n:\Sigma^{-f(n)}BP_\ast/I_n \longr BP_\ast/I_n^\infty
$$
where $f(1) = 0$ and for $n > 1$
$$
f(n) = 2(p-1) + 2(p^2-1) + \cdots +2(p^{n-1}-1).
$$
The morphism $\iota_n$ sends the generator of $BP_\ast/I_n$ to the
element
$$
\frac{1}{pv_1\cdots v_{n-1}} \in BP_\ast/I_n^\infty.
$$
This can also be defined inductively by letting $\iota_0:BP_\ast \to
p^{-1}BP_\ast$ be the usual inclusion and noting there is a commutative
diagram of comodules
$$
\xymatrix{
\Sigma^{-f(n)}BP_\ast/I_{n} \rto^{v_n} \dto_{\iota_{n}} &
\Sigma^{-f(n+1)} BP_\ast/I_n \dto\\
BP_\ast/I_{n}^\infty \rto & v_{n}^{-1}BP_\ast/I_n^\infty
}
$$
and taking the induced map on quotients. The element
$v_n \in BP_\ast/I_n$ is primitive;
this is a consequence of the formula \ref{p-series-n} and the fact
that if a formal group law $F$ has height at least $n$, then $v_n^F$ 
is a strict isomorphism invariant. This yields a sub-ring
$$
\FF_p[v_n] \subseteq H^{0,\ast}(BP_\ast/I_n).
$$
Since $H^{0,\ast}$ is left exact (preserves injections) we get non-zero
classes $v_n^k$ in
$$
H^{0,2k(p^n-1) - f(n)}(BP_\ast/I_n^\infty) \subseteq
H^{0,2k(p^n-1) - f(n)}(v_n^{-1}BP_\ast/I_n^\infty).
$$
The image of the class $v_n^k$ in 
$$
H^{0,2k(p^n-1) - f(n)}(BP_\ast) =
\Ext^s_{BP_\ast BP}(\Sigma^{2k(p^n-1) - f(n)}BP_\ast,BP_\ast)
$$
is the $k$th element in $n$th Greek letter family and is sometime
written $\alpha_k^{(n)}$. There is no reason at this point to suppose
that these classes are non-zero -- although we will see shortly that
for low values of $n$ there can be no more differentials in the
chromatic spectral sequence.

If $n$ is small we use the actual
Greek letter. Thus the cases $n=1$ and $n=2$ are written
$$
\alpha_k \in \Ext^1_{BP_\ast BP}(\Sigma^{2k(p-1)}BP_\ast,BP_\ast)
$$
and
$$
\beta_k \in \Ext^2_{BP_\ast BP}(\Sigma^{2k(p^2-1) - 2(p-1)}BP_\ast,BP_\ast).
$$
\end{exam}

This brings us to our first periodic families in the stable homotopy
groups of spheres.

\begin{thm}\label{period1}Let $p\geq 3$. Then the class
$$
\alpha_k \in \Ext^1_{BP_\ast BP}(\Sigma^{2k(p-1)}BP_\ast,BP_\ast)
$$
is a non-zero permanent cycle in the ANSS spectral sequence
and detects an element of order $p$. 

Let $p\geq 5$. Then the class
$$
\beta_k \in \Ext^2_{BP_\ast BP}(\Sigma^{2k(p^2-1) - 2(p-1)}BP_\ast,BP_\ast).
$$
is a non-zero permanent cycle in the ANSS spectral sequence
and detects an element of order $p$. 
\end{thm}.

At the prime $2$ we know $d_3\alpha_3 \ne 0$; at the prime $3$
we know $d_5\beta_4 \ne 0$.

We will see in the next section that these classes are non-zero 
at $E_2$. They cannot be boundaries, but to show they are permanent
cycles requires some homotopy theory. 
\begin{rem}[{\bf Smith-Toda complexes}]\label{smith-toda}Suppose
there exists finite $CW$ spectrum $V(n)$ so that there is an 
isomorphism of comodules
$$
BP_\ast V(n) \cong BP_\ast/I_n
$$
and so that there is a self map $v:\Sigma^{2(p^n-1)}V(n) \to V(n)$
which induces multiplication by $v_n$ on $BP_\ast(-)$. A short
calculation shows that the top cell of $V(n)$ must be in degree
$f(n) + n$ -- the degree of the Milnor operation $Q_0\ldots Q_{n-1}$.
Let $i:S^0 \to V(n)$ and $q:V(n) \to S^{f(n)+n}$ be the inclusion
of the bottom cell and the projection onto the top cell respectively.
Then the composition
$$
\xymatrix{
S^{2k(p^n-1)} \rto^-i & \Sigma^{2k(p^n-1)}V(n) \rto^-{v^k} &
V(n) \rto^-q & S^{f(n)+n}
}
$$
is detected by $\alpha_k^{(n)}$ in the ANSS; thus, this class is
a permanent cycle and, if non-zero, detects a non-trivial element
in the homotopy groups of spheres. These ideas originated with
\cite{lsmith} and \cite{toda} and have been explored in great
details by many people; see \cite{Rav} \S 5.5.

This pleasant story is marred by the fact that the $V(n)$ are known
to exist only for small values of $n$ and large values of $p$.
The spectrum $V(0)$ is the Moore spectrum, which exists for 
all $n$, but has a $v_1$-self map only if $p > 2$. (At $p=2$
one can realize $v_1^4$, which is of degree $8$, a number which
should not surprise Bott periodicity fans.) The spectrum $V(1)$
exists if $p > 2$, but has a $v_2$-self map only if $p>3$. One
can continue for a bit, but the spectrum $V(4)$ is not known to
exist at any prime and non-existence results abound. See
especially \cite{nave}.

Note that any $v_n$-self map must have the property that $v^k \ne 0$
for all $k$. The existence and uniqueness of maps of finite
complexes which realize multiplication by powers of $v_n$ is
part of the whole circle of ideas in the nilpotence conjectures.
See \cite{nilpotence}, \cite{DHS}, \cite{HS}, and \cite{orange}.
\end{rem}

\begin{rem}\label{monochrom}The groups
$$
\Ext_{BP_\ast BP}^s(\Sigma^t BP_\ast,v_n^{-1}BP_\ast/I_n^\infty)
$$
form the $E_2$-term of an ANSS converging to an appropriate
suspension of the fiber of the map between localizations
$$
L_n S^0 \longr L_{n-1} S^0.
$$
This fiber, often written $M_nS^0$ is the ``monochromatic'' piece
of the homotopy theory of the sphere $S^0$. Again see
\cite{orange}.
\end{rem}

 \section{Change of rings}

We now begin trying to calculate $E_2$-terms of the ANSS spectral
sequence based on some complex oriented, Landweber exact
homology theory $E_\ast$ over $\ZZP$. The primordial examples
are $BP_\ast$ and $\ZZP \otimes MU_\ast$; the Cartier idempotent
shows that the Hopf algebroids
$$
(BP_\ast,BP_\ast BP)\qquad\mathrm{and}\qquad 
\ZZP \otimes (MU_\ast,MU_\ast MU)
$$
represent equivalent affine groupoid schemes over $\ZZP$, from which
is follows that
$$
\ZZP \otimes \Ext^s_{MU_\ast MU}(\Sigma^t MU_\ast, MU_\ast)
\cong \Ext^s_{BP_\ast BP}(\Sigma^t BP_\ast, BP_\ast ).
$$
This is an example of a {\it change of rings} theorem; we give 
several such theorems in this section.

Recall that a 
functor $f:G \to H$ of groupoids is an equivalence if
it induces a weak equivalence on classifying spaces $BG \to BH$;
equivalently, it induces an isomorphism
$$
f_\ast: \pi_0 G \mathop{\longr}^{\cong} \pi_0 H
$$
on the sets of isomorphism classes of objects and for each
object $x \in G$ we get an isomorphism
$$
f_\ast: \pi_1(G,x) = \Aut_G(x) \mathop{\longr}^{\cong} \Aut_H(f(x))
= \pi_1(H,f(x)).
$$
This notion has a weakening. 

\begin{defn}\label{local-fl-top} Suppose $G$ and $H$ are actually
functors to groupoids from commutative rings. Then a functor
$f:G \to H$ (now a natural transformation) is an {\bf equivalence
locally in the flat topology} if 
\begin{enumerate}

\item for all rings $R$ and all $x \in H(R)$ there exists a faithfully
flat extension $R \to S$ of rings and $y \in G(R)$ so that
$f(y) \cong x$ in $H(S)$;

\item for all rings $R$ and all $x,y \in G(R)$ and $\phi:f(x) \to f(y)$
in $H(R)$, there is a faithfully flat extension $R \to S$ and
a unique $\psi:x \to y$ in $G(S)$ so that $f\psi = \phi$ in $H(S)$.
\end{enumerate}
\end{defn} 

It is equivalent to say that the induced morphism $f_\ast:\pi_iH \to \pi_iG$
($i=0,1$)
of presheaves on affine schemes induces an isomorphism on
associated sheaves in the flat topology.

\begin{thm}[{\bf Change of Rings I}]\label{change-of-rings}Let $(A,\Gamma) \to (B,\La)$
be a morphism of Hopf algebroids. If $M$ is an $(A,\Gamma)$-comodule,
then $B\otimes_A M$ is a $(B,\La)$-comodule and there is
a homomorphism 
$$
\Ext_\Gamma^s(\Sigma^tA,M) \longr \Ext^s_\La(\Sigma^tB,B \otimes_AM).
$$
If the morphism of Hopf algebroids induces an equivalence of affine
groupoid schemes locally in the flat topology, then this morphism
on $\Ext$-groups is an isomorphism.
\end{thm}

\begin{rem}\label{change-of-rings-proof}
A conceptual proof of this can be constructed using the ideas
of Remark \ref{qc}; indeed, we can even relax the requirement
that we are working with Hopf algebroids and consider groupoid
schemes such as $\cU(n)$. See \cite{groupoid}.
\end{rem}

\begin{exam}\label{flat-base-ch}For example, let $(A,\Ga)$ be
a Hopf algebroid and $A \to B$ a faithfully flat extension of rings.
Then we get a new Hopf algebroid $(B,\La)$ with
$$
\La = B \otimes_A \Ga \otimes_A B
$$
and we have that for all $(A,\Ga)$-comodules $M$
$$
\Ext^s_\Ga(\Sigma^tA,M) \cong \Ext^s_{\La}(\Sigma^t B,B \otimes_A M).
$$
To be very concrete let $(A,\Ga)=(L,W)=(MU_\ast,MU_\ast MU)$, the Hopf ring representing
the groupoid scheme of graded formal group laws.
Let $B = L[u^{\pm 1}]$ where $u$ is a unit in degree $2$. Then
$$
B = L_0[u^{\pm 1}]
$$
where $L_0$ is now the Lazard ring concentrated in degree $0$.
We also have
$$
B \otimes_L W \otimes_L B= L_0[u^{\pm 1},b_0^{\pm1 },b_1,\ldots] =
W_0[u^{\pm 1}]
$$
where the $b_i$ are the coefficients of the universal isomorphism
of (ungraded) formal group laws. This yields
$$
\Ext^s_W(\Sigma^tL,M) \cong
\Ext^s_{W_0[u^{\pm 1}]}(\Sigma^t L_0[u^{\pm 1}],M[u^{\pm 1}]).
$$
This can be further simplified using Example \ref{2-period}.
\end{exam}

\begin{exam}[{\bf Morava $K$-theory to group cohomology}]\label{grp-coh}
We apply these results to the Hopf algebroid 
$(K(n)_\ast,\Sigma(n))$ constructed in Example \ref{Morava-stabilizer}.
Let 
$$
K(n)_\ast = \FF_p[v_n^{\pm 1}] \longr \FF_{p^n}[u^{\pm 1}]
$$
be the faithfully flat extension with $u^{p^n-1} = v_n$. The
extended Hopf algebroid is
$$
\FF_{p^n}[u^{\pm 1}] \otimes_{\FF_{p^n}} (\FF_{p^n},
\FF_{p^n} \otimes S(n) \otimes \FF_{p^n}).
$$
(We must take a little care with this formula, as $\eta_R(u) = b_0u$,
where $b_0$ is the leading term of the universal isomorphism
of ungraded formal groups laws.)
By Galois theory, we know
$$
\FF_{p^n} \otimes \FF_{p^n} \cong \map(\Gal(\FF_{p^n}/\FF_p),\FF_{p^n})
$$
and
$$
\FF_{p^n} \otimes S(n) \cong \map(\Aut(\Ga_n),\FF_{p^n})
$$
by Example \ref{Morava-stabilizer}. The group $\Aut(\Ga_n) = S_n$
is a profinite group and the last set of maps is actually the continuous
maps. Putting these two facts together we have
$$
\FF_{p^n} \otimes S(n) \otimes \FF_{p^n} \cong
\map(G_n,\FF_{p^n})
$$
where $G_n = \Gal(\FF_{p^n}/\FF_p) \rtimes S_n$ is the semi-direct
product. It follows that if $M$ is a $(K(n)_\ast,\Sigma(n))$
comodule concentrated in even degrees (see Example \ref{2-period})
then
\begin{align*}
\Ext_{K(n)_\ast \otimes S(n)}(\Sigma^{2t}K(n)_\ast,M)
&\cong H^s(G_n,\FF_{p^n} \otimes M_0(t))\\
&\cong H^s(S_n,\FF_{p^n} \otimes M_0(t))^{\Gal(\FF_{p^n}/\FF_n)}.
\end{align*}
Cohomology here is continuous cohomology. If $N$ is continuous
$G_n$ module, $\phi \in S_n$ and $u^t x \in N(t)$, then
$$
\phi(u^tx) = \phi'(0)^tu^t (\phi x).
$$
\end{exam}

\def\leten{{{E_n}}}
\def\reten{{{R(\FF_{p^n},\Ga_n)}}}

\begin{exam}\label{grp-coh-2}The previous example can be extended
to a wider class of examples. Let $M$ be a comodule over
$(E(n)_\ast,E(n)_\ast E(n))$ and suppose $n$ is $I_n$ torsion.
Here $I_n = (p,v_1,\cdots,v_{n-1}) \subseteq E(n)_\ast$ and
$$
E(n)_\ast/I_n \cong K(n)_\ast.
$$
We'd like to relate the Ext groups of $M$ to an appropriate
group cohomology. Let $(E_n)_\ast = \leten_\ast$ be the $2$-periodic
Lubin-Tate homology theory obtained from the deformation theory
of the formal group law $\Ga_n$ over $\FF_{p^n}$. Thus 
$\leten_\ast = \reten[u^{\pm 1}]$ where $u$ is
in degree $2$ and the deformation ring
$$
\reten = W(\FF_{p^n})[[u_1,\cdots,u_{n-1}]]
$$
is in degree $0$. This is a local ring with maximal ideal
$$
\mm = (p,u_1,\ldots,u_{n-1}).
$$
Define a faithfully flat extension $E(n)_\ast \to \leten_\ast$ by
$$
v_i \mapsto \brackets{u_iu^{p^i-1},}{i < n;}{u^{p^n-1},}{n=1.}
$$
Then $\mm = [\leten_\ast \otimes_{E(n)_\ast} I_n]_0$. Let's also
write $\leten_\ast\leten$ for the completion of
$$
\leten_\ast \otimes_{E(n)_\ast} E(n)_\ast E(n) \otimes_{E(n)_\ast} \leten_\ast
$$
at the maximal ideal $\mm$. Then Lubin-Tate theory tells us that
we obtain the continuous maps:
$$
\leten_\ast\leten \cong \map(G_n,\leten_\ast).
$$
Let $M$ be an $(E(n)_\ast,E(n)_\ast E(n))$-comodule which is $I_n$-torsion
as a $E(n)_\ast$-module.
Set $\bar{M} = \leten_\ast \otimes_{E(n)_\ast} M$.
Then, because of the torsion condition,
\begin{align*}
 \map(G_n,\bar{M}) &\cong \leten_\ast \leten \otimes_{\leten_\ast} \bar{M}\\
&\cong 
 \leten_\ast \otimes_{E(n)_\ast} E(n)_\ast E(n) \otimes_{E(n)_\ast} \leten_\ast
\otimes_{\leten_\ast} \bar{M}
\end{align*}
Here $\map(G_n,\bar{M})$ is the module of continuous maps, where
$\bar{M}$ has the discrete topology. From this, and argument
similar to that of the previous example, we have
\begin{equation}\label{reduce-to-coh-2}
\Ext^s_{E(n)_\ast E(n)}(\Sigma^{2t},M) \cong
H^s(G_n,\bar{M}_0(t)).
\end{equation}
\end{exam}

The final change of rings result is the following -- this is Morava's
original change of rings result. See \cite{devinatz}, \cite{Morava} and \cite{Rav}.

\begin{thm}[{\bf Change of rings II}]\label{change-2}Let
$M$ be a $(BP_\ast,BP_\ast BP)$-comodule which is
$I_n$-torsion as $BP_\ast$-module. Let
$$
\bar{M} = \leten_\ast \otimes_{BP_\ast} M
\cong  \leten_\ast \otimes_{E(n)_\ast} E(n)_\ast \otimes_{BP_\ast} M.
$$
Then if $M$ is concentrated in degrees congruent to $0$ modulo
$2(p-1)$, there is a natural isomorphism
$$
\Ext_{BP_\ast BP}^s(\Sigma^{2t}BP_\ast,v_n^{-1}M)
\cong H^s(G_n,\bar{M}_0(t)).
$$
The odd cohomology groups are zero.
\end{thm}

\begin{rem}\label{change-2-proof}Here is a sketch of the proof. It
is sufficient, by Example \ref{grp-coh-2} to show
$$
\Ext_{BP_\ast BP}^s(\Sigma^{2t}BP_\ast,v_n^{-1}M)
\cong \Ext_{E(n)_\ast E(n)}(\Sigma^{2t}E(n)_\ast, E(n)_\ast
\otimes_{BP_\ast} M).
$$
We use the ideas of Remark \ref{qc}. Let $i:\cU(n) \to \fglp$
be the inclusion of the open sub-groupoid scheme of
formal group laws of height no more than $n$. The pull-back
functor $i^\ast$ from quasi-coherent sheaves (that is, comodules)
over $\fglp$ has a left exact right adjoint $i_\ast$. Then,
by Remarks \ref{filt-basic}.3 and \ref{change-of-rings-proof}, 
it is equivalent to prove that
$$
H^\ast(\fglp,v_n^{-1}\cF_M) \cong H^\ast (\cU(n),i^\ast \cF_M).
$$
I have written $\cF_M$ for the quasi-coherent sheaf
associated to $M$ and suppressed the internal grading. 
It turns out that because $M$ is
$I_n$-torsion
$$
i_\ast i^\ast \cF_M \cong v_n^{-1}\cF_M.
$$
Indeed, this formula is an alternate and conceptual proof
of Lemma \ref{invert-vn}. Furthermore, 
$$
R^qi_\ast i^\ast \cF_M = 0,\qquad q > 0.
$$
The result now follows from the composite functor spectral
sequence
$$
H^p(\fglp,R^qi_\ast \cE) \Longrightarrow H^{p+q}(\cU(n),\cE).
$$
\end{rem}

\section{The Morava stabilizer group}\label{sect-group}

The Morava stabilizer group $G_n$ was introduced in Example
\ref{Morava-stabilizer} and now we write down a little bit
about its structure. Much of this material has classical
roots; for more discussion see \cite{Rav} \S A.2.2 and
\cite{hwhduke}.

We have fixed an ungraded $p$-typical formal
group law $\Ga_n$ over $\FF_p$ with $p$-series
$$
[p](x) = x^{p^n}.
$$
(The choice of $\Ga_n$ is simply to be definite; there is
a general theory which we will not invoke.) Then $G_n$
is the set of pairs $(\sigma,\phi)$ where
$\sigma:\FF_{p^n} \to \FF_{p^n}$ is a element of the
Galois group of $\FF_{p^n}$ over $\FF_p$ and
$$
\phi:\Ga_n \longr \sigma^\ast \Ga_n
$$
is an non-strict isomorphism of formal groups. We see immediately
that we have a semi-direct product decomposition
$$
G_n \cong \Gal(\FF_{p^n}/\FF_p) \rtimes \Aut(\Ga_n)
$$
where we are write $\Aut(\Ga_n)$ for the automorphisms of
$\Ga_n$ over $\FF_{p^n}$. These automorphisms are
the units in the ring $\End(\Ga_n)$ of endomorphisms
and any such endomorphism $f(x)$ can be written
$$
f(x) = a_0 x +_{\Ga_n} a_1x^{p} +_{\Ga_n} a_2 x^{p^2} +_{\Ga_n} \cdots
$$
where $a_i \in \FF_{p^2}$. The automorphisms have $a_0 \ne 0$.

Of particular note are the endomorphisms $p=[p](x) = x^{p^n}$ and
$S(x) = x^p$. We have, of course, that $S^n = p$.

The subring of the endomorphisms with elements
$$
f(x) = a_0 x +_{\Ga_n} a_{n}x^{p^n} +_{\Ga_n} a_{2n}
 x^{p^{2n}} +_{\Ga_n} \cdots
$$
is a commutative complete local ring with with maximal ideal
generated by $p$ and residue field $\FF_{p^n}$; thus, from the
universal property of Witt vectors, it is canonically isomorphic
to $W(\FF_{p^n})$. It follows that
$$
\End(\Ga_n)\cong W(\FF_{p^2})\langle S \rangle/(S^n-p).
$$
The brackets mean that $S$ does not commute the Witt vectors;
in fact $Sa = \phi(a)S$ where $\phi$ is the Frobenius. While we
will not use this fact, the endomorphism ring has a classical
description as the maximal order in a division algebra over
the $p$-adic rationals. We now have 
$$
\Aut(\Ga_n) \cong [W(\FF_{p^2})\langle S \rangle/(S^n-p)]^\times.
$$
Thus an element $\Aut(\Ga_n)$ can be expressed as a sum
$$
b_0 + b_1S + \cdots b_{n-1}S^{n-1}
$$
where $b_i \in W(\FF_{p^n})$ and $b_i \ne 0$ modulo $p$.

\begin{rem}\label{morava-basics}
Here are a few basic facts about this group.
\begin{enumerate}

\item The powers of $S$ define a filtration on $\End(\Ga_n)$
and $\End(\Ga_n)$ is complete with respect to this filtration.
There is an induced filtration on $\Aut(\Ga_n)$ with
$$
F_{i/n}\Aut(\Ga_n) = \{\ f\ |\ f \equiv 1\ \mathrm{mod}\ S^{i}\ \}.
$$
In this filtration the group is a profinite group. The fractional
notation is so that $p$ has filtration degree $1$. 

\item We define $S_n = F_{1/n}\Aut(\Ga_n)$; it is the set of
automorphisms $f(x)$ with $f'(0) = 1$ -- the strict automorphisms.
This is the $p$-Sylow subgroup (in the profinite sense) of
the automorphism group.

\item The elements $a_0x \in \Aut(\Ga_n)$ define a subgroup
isomorphic to $\FF_{p^2}^\times$ splitting the projection
\begin{align*}
\Aut(\FF_{p^n}) &\longr \FF_{p^2}^\times\\
f &\mapsto f'(0).
\end{align*}
Thus we have a
semi-direct product decomposition
$$
\Aut(\Ga_n) \cong \FF_{p^n}\rtimes S_n.
$$

\item The center of $C \subseteq S_n$ is the subgroup of 
elements of the form
$b_0$ where $b_0 \in \ZZ_p^\times \subseteq W(\FF_{p^n})^\times$
and $b_0 \equiv 1$ modulo $p$. As power series these are 
of the form
$$
f(x) =  x +_{\Ga_n} a_{n}x^{p^n} +_{\Ga_n} a_{2n}
 x^{p^{2n}} +_{\Ga_n} \cdots
$$
where $a_i \in \FF_p$. Thus $\ZZ_p \cong C \subseteq \ZZ_p^\times$.

\item The action of $\Aut(\Ga_n)$ on $\End(\Ga_n)$ defines a 
homomorphism
$$
\Aut(\Ga_n) \to \mathrm{Gl}_n(W(\FF_{p^n}));
$$
taking the determinant defines a homomorphism
$$
S_n \longr \ZZ_p^\times.
$$
(A priori, the determinant lands in $W(\FF_{p^n})^\times$, but
we check it actually lands in the subgroup $\ZZ_p^\times$.) 
The elements of the subgroup $S_n$ will all map
into $\ZZ_p$ and the composition $\ZZ_p \cong C \to S_n \to \ZZ_p$ is
multiplication by $n$; hence, if $(n,p) = 1$, we get a splitting
$$
S_n \cong C \times S_n^1
$$
where $S_n^1$ is the kernel of determinant map $\Aut(\Ga_n)
\to \ZZ_p^\times$.

\item A deeper part of the theory studies the finite subgroups of
$S_n$. We note here only that $S_n$ contains elements of
order $p$ if and only if $p-1$ divides $n$. At the prime $2$
this is bad news, of course.
\end{enumerate}
\end{rem}

We now turn to a discussion of the continuous cohomology
of these groups. The first statement of the following
-- and much, much more -- can be found in \cite{lazgroup}. The
second is follows from the last item of the previous remark.
If $G$ a $p$-profinite group is a Poincar\'e duality group
of dimension $m$, then $H^\ast(G,M)$ (with any
coefficients) is a Poincar\'e duality algebra with top
class in degree $m$.

\begin{prop}\label{vfcd}The group
$S_n$ has a finite index subgroups $H_n$ which is  a
Poincar\'e duality group of dimension $n^2$. If
$p-1$ does not divide $n$, we may take $H_n$ to
be the entire group.
\end{prop}

Let's now get concrete and try to calculate
$H^\ast (G_2,\FF_{p^2}(\ast))$, $p > 2$.
The action of $G_2$ on $\FF_{p^2}(\ast)$ is through the quotient
subgroup
$$
\Gal(\FF_{p^2}/\FF_p) \rtimes \FF_{p^2}^\times.
$$
The kernel of quotient map is the $p$-Sylow subgroup $S_2$,
which (since $p>2$) can be written as
a product
$$
C \times S_2^1$$
where $\ZZ_p\cong C \subset S_2$ is the center. 

If $G$ is a profinite group and $\FF$ is a field on which
$G$ acts trivially, $H^1(G,\FF)$ is the group of $1$-cocyles;
furthermore,
these are simply the continuous group homomorphisms
$$
f:G \longr \FF.
$$

Define a $1$-cocyle $\zeta:C \to \FF_p$ on the
center by
$$
\zeta(f) = a_2.
$$
Here we are using the notation of Remark \ref{morava-basics}.4.

\begin{lem}\label{center}There natural map
$$
E(\zeta) \longr H^\ast(C,\FF_{p^2})
$$
is an isomorphism. Furthermore as an $\FF_{p^2}^\times$-module
$$
H^1(C,\FF_{p^2}) \cong \FF_{p^2}(0).
$$
\end{lem}

\begin{proof}The first part is standard; the second part follows from
the fact that $C$ is the center, so conjugation by an element
$a_0x$, $a_0 \in \FF_{p^2}$ is trivial.
\end{proof}

Now define a 1-cocycles $\tau_1,\tau_2:S^1_2 \to \FF_{p^2}$
on $S_2$ by
$$
\tau(f) = a_1\qquad\mathrm{and}\qquad \tau_2(f) = a_1^p.
$$
Again again we use the notation of  Remark \ref{morava-basics}. 

\begin{lem}\label{1cohs21}There is an isomorphism
of $\FF_{p^2}^\times$-modules
$$
H^1(S_2^1,\FF_{p^2}) \cong \FF_{p^2}(p-1) \oplus \FF_{p^2}(p^2-p)
$$
with summands generated by $\tau_1$ and $\tau_2$ respectively.
Furthermore
$$
\tau_1\tau_2 = 0.
$$
\end{lem}

\begin{proof} We check that 
$$
\tau_1:S_2^1/[S_2^1,S_2^1] \cong \FF_{p^2}
$$
is an isomorphism. Then
$$
H^1(S_2^1,\FF_{p^2}) \cong \Hom(S_2^1/[S_2^1,S_2^1],\FF_{p^2})
\cong \FF_{p^2} \times \FF_{p^2}
$$
generated by $\tau_1$ and the Frobenius applied to $\tau_1$ --
that is, $\tau_2$. To get the action of $\FF_{p^2}$, let
$a_0x$ be typical element of that group and $f(x) \in S_2^1$.
Then we check
$$
\tau_1(a_0^{-1}(f(a_0x)) = a_0^{p-1}\tau_1(f(x)),
$$
from which it follows that
$$
\tau_2(a_0^{-1}(f(a_0x)) = a_0^{p^2-p}\tau_2(f(x)),
$$
as claimed. 

To show the product is zero we define $\sigma(f(x)) = a_2$.
(This is not the same of $\zeta$, even though it has the same formula,
as the domain of definition is $S_2^1$, not the center.)
Then there is an equality of $2$-cocyles $\partial\sigma = -\tau_1\tau_2$.
\end{proof}

We can now define elements $\gamma_1$ and $\gamma_2$
in $H^2(S_2^1,\FF_{p^2})$ by
\begin{align*}
\gamma_1 &= \langle \tau_1,\tau_2,\tau_1 \rangle\\
\gamma_2 &= \langle \tau_1,\tau_2,\tau_1 \rangle
\end{align*}
This next result can be
proved using Hopf algebra techniques as in \cite{Rav} Theorem 6.3.22
or using  Lazard's work \cite{lazgroup}. See \cite{hwhduke} for
an introduction to the latter approach. We must
restrict to primes $p > 3$ as $G_2$ at the prime $3$ has
$3$-torsion.

\begin{prop}\label{hs21}Let $p > 3$. The cohomology ring $H^\ast(S_2^1,\FF_{p^2})$
has generators $\tau_1$,$\tau_2$,$\gamma_1$, and $\gamma_2$ with
all products zero except
$$
\tau_1\gamma_2 = \gamma_1\tau_2
$$
generates $H^3(S_2^1,\FF_{p^2})$.
\end{prop}

The next result is a matter of taking invariants and Lemma \ref{1cohs21}.
Note that
$$
H^0(G_2,\FF_{p^2}(\ast)) = H^0(\Gal(\FF_{p^2}/\FF_p \rtimes \FF_{p^2}^\times,
\FF_{p^2}[u^{\pm 1}]) \cong \FF_p[v_2^{\pm 1}]
$$
where $v_2 = u^{p^2-1}$.

\begin{cor}\label{1cohG2}If $p \geq 3$, the free $\FF_p[v_2^{\pm 1}]$-module
$$
H^1(G_2,\FF_{p^2}(\ast))
$$
is generated by elements
\begin{align*}
\zeta &\in H^1(G_2,\FF_{p^2}(0))\\
h_0 &\in  H^1(G_2,\FF_{p^2}(p-1))\\
h_1 &\in H^1(G_2,\FF_{p^2}(p^2-p))\\
\end{align*}
\end{cor}

\section{Deeper periodic phenomena}

In this last section we calculate (or at least talk about calculating)
enough of the chromatic spectral sequence \ref{chrom-ss} in
order to understand the $0$, $1$, and $2$-lines of the
ANSS. At larger primes (here $p > 3$), this is largely an
algebraic calculation. We begin with the easy calculation:

\begin{lem}\label{0-line}For all primes $p$,
$$
H^{s,t}(p^{-1}BP_\ast) \cong \brackets{\QQ,}{t=0;}{0}{t\ne 0}.
$$
\end{lem}

\begin{proof}This is a change of rings arguments. We check
that
$$
H^{\ast,\ast}(p^{-1}BP_\ast) \cong
\Ext_{\QQ \otimes W}^\ast(\Sigma^\ast(\QQ \otimes L),\QQ \otimes L)
$$
and then apply the change of rings theorem \ref{change-of-rings}
to the morphism of Hopf algebroids
$$
\QQ \otimes (L,W) \longr (\QQ,\QQ \otimes_L W \otimes_L \QQ)
$$
induced by the map $L \to \QQ$ classifying the additive formal group
law. Since the additive formal group law over the rational numbers has
no non-trivial strict isomorphisms we deduce that
$$
\QQ \otimes_L W \otimes_L \QQ \cong \QQ.
$$
The result follows.
\end{proof}

For the higher entries in the chromatic spectral sequence we use
Morava's change of ring theorem \ref{change-2}. Let
$$
R_n = \reten = (E_n)_0 \cong W(\FF_{p^n})[[u_1,\cdots,u_{n-1}]]
$$
be the Lubin-Tate deformation ring and let 
$$
\mm = (p,u_1,\cdots,u_{n-1}) \subseteq \reten
$$
be its maximal ideal. Then we have
\begin{equation}\label{the-big-iso}
H^{s,2t}(v_n^{-1}BP/I_n^\infty) \cong H^s(G_n,(R_n/\mm^\infty)(t)).
\end{equation}

We now look at $n=1$. The ring $R_1$ is isomorphic to the
$p$-adic integers, $\mm = (p)$ and $R_1/m^\infty \cong \ZZ/p^\infty$.
The group $G_1 = \ZZ_p^\times$, the units in the $p$-adic
integers,
and the action $G_1$ on $\ZZ/p^\infty(t)$ is given by
$$
a(u^tx) = u^t(a^tx).
$$
In particular the action on $\ZZ/p^\infty(0)$ is trivial.
The calculation for the $1$ line begins with the following result.

\begin{thm}\label{1-line}Let $p > 2$. Then $H^{s,2t} =
H^s(G_1,\ZZ/p^\infty(t))=0$ unless $s=0$ or $1$ and
$2t=2k(p-1)$ for some $k$. In the case $k=0$ we have
$$
H^{0}(G_1,\ZZ/p^\infty) = H^1(G_1,\ZZ/p^\infty) \cong \ZZ/p^\infty
$$
and if $k = p^rk_0$ with $(p,k_0)=1$,
$$
H^0(G_1,\ZZ/p^\infty(t)) \cong \ZZ/p^{r+1}\ZZ
$$
and $H^1(G_1,\ZZ/p^\infty(t))=0$.
\end{thm}

\begin{proof}We can write $\ZZ_p^\times \cong C_{p-1} \times \ZZ_p$
where $C_{p-1} \cong \FF_p^\times$ is cyclic of order $p-1$. Then
\begin{equation}\label{red-to-z}
H^s(G_1,\ZZ/p^\infty(t)) \cong H^s(\ZZ_p,\ZZ/p^\infty(t))^{C_{p-1}}
\end{equation}
is zero unless $t = k(p-1)$; in this case,
$$
H^s(\ZZ_p,\ZZ/p^\infty(t))^{C_{p-1}} \cong H^s(\ZZ_p,\ZZ/p^\infty(t)).
$$
Choose a generator $\gamma \in \ZZ_p \subseteq \ZZ_p^\times$; since
$p$ is odd, this
is a $p$-adic unit $x$ congruent to $1$ modulo $p$, and so that  $x-1$
is non-zero modulo $p^2$. Then the cohomology
can be calculated using the very short cochain complex
$$
\gamma^t -1:\ZZ/p^\infty(t) \longr \ZZ/p^\infty(t).
$$
\end{proof}

\begin{rem}\label{final-zero}The morphism
$$
\QQ \cong H^{0,0}(p^{-1}BP_\ast) \to H^{0,0}(v_1^{-1}BP_\ast) \cong
\ZZ/p^\infty
$$
is surjective as the target group is generated by the elements
$1/p^t \in [v_1^{-1}BP_\ast]_0$. Thus
$$
\pi_0S^0 \cong \Ext_{BP_\ast BP}^0(BP_\ast,BP_\ast) \cong \ZZP,
$$
as it better.
\end{rem}

\begin{rem}\label{one-at-two}At the prime $2$, we have an 
isomorphism
$$
G_1 \cong \ZZ_2^\times \cong C_2 \times \ZZ_2
$$
where $C_2 = \{\pm 1\}$ and $\ZZ_2 \cong \ZZ_2^\times$
are the $2$-adic units with are congruent to $1$ modulo $4$.
The decomposition of \ref{red-to-z} now becomes a spectral
sequence, which collapses but yields $H^{s,\ast}(v_1^{-1}BP_\ast)
\ne 0$ for all $s$.
\end{rem}

\begin{rem}[{\bf The image of $J$}]\label{final-one} Let $p>2$.
The elements in non-positive degree
in $H^{0,\ast}(v_1^{-1}BP_\ast)$ cannot be permanent cycles in the
chromatic spectral sequence as we know $H^{1,\ast}(BP_\ast)=0$
in those degrees -- see Lemma \ref{first-conv}. The elements
of positive degree are permanent cycles, indeed if $2t = 2k(p-1)$
and $k = p^rk_0$, then
$$
v_1^k/p^{r+1} \in H^{0,2k(p-1)}(BP_\ast/(p^\infty))
$$
maps to a generator of $H^{0,2k(p-1)}(v_1^{-1}BP_\ast/(p^\infty))$.
Here we use that
$$
\eta_R(v_1) = v_1 + px
$$
to get $v_1^k/p^{r+1}$ to be primitive. Thus we have an isomorphism
$$
H^{1,2k(p-1)}(BP_\ast) \cong \ZZ/p^{r+1}\ZZ
$$
generated by an element $\alpha_{k/r}$. These elements
are all permanent cycles on the ANSS, as well. This can
be seen either by using a variant of the Smith-Toda argument
of Remark \ref{smith-toda} with the mod $p^{r+1}$ Moore spectra
$S^0 \cup_{p^{r+1}} D^1$ or noting that these elements
detect the image of the $J$ homomorphism. 
\end{rem}

\begin{rem}[{\bf Hopf invariant one}]\label{hopf} Let $p > 2$. The evident
map of ring spectra $BP \to H\ZZ/p$ induces a commutative
diagram of spectral sequences
$$
\xymatrix{
\Ext^s_{BP_\ast BP}(\Sigma^t BP_\ast,BP) \ar@2{->}[r]\dto
& \pi_{t-s}S^0_{(p)}\dto\\
\Ext^s_{A_\ast}(\Sigma^t \FF_p,\FF_p) \ar@2{->}[r]& \pi_{t-s}(S^0)^\cmpl_p.
}
$$
The map
$$
\Ext^1_{BP_\ast BP}(\Sigma^t BP_\ast,BP)\longr
\Ext^1_{A_\ast}(\Sigma^t \FF_p,\FF_p)
$$
is almost entirely zero: $\alpha_{k/r}$ is in the kernel unless
$k=r=1$. As a result, the Hopf invariant $1$ elements
$$
h_i \in \Ext^1_{A_\ast}(\Sigma^{p^i} \FF_p,\FF_p)
$$
must support differentials if $i > 0$ -- if they were permanent
cycles there would detect a class of filtration $0$ or $1$ in
the ANSS. A similar argument works at the prime $2$.
Since this argument uses only information from chromatic
level $1$, it could be done just as well with $K$-theory; therefore,
this is a reformulation of Atiyah's proof of the non-existence
of elements of Hopf invariant $1$.

While I've not yet talked about the $2$-line, there is
a very large kernel at that line as well. Indeed, using the generators
of \ref{beta-names} below, only $\beta_2$, $\beta_{p^i/p^i}$, and
$\beta_{p^i/p^i-1}$ do not map to zero.
\end{rem}

\begin{rem}\label{p-adic-nghbd}In Remark \ref{final-one} we wrote down
specific generators $v_1^k/p^{r+1}$ (with $k = p^rk_0$) for
$H^{0,\ast}(BP_\ast/(p^\infty))$. These have the property
that as $r \to \infty$, they generate a torsion group
of higher and higher order and approach the case of $k=0$,
which yields the infinite divisible group $\ZZ/p^\infty$. The
requirement that $r \to \infty$ is the same as the requirement
that $k \to 0$ in the $p$-adic topology on the integers. Thus
we regard the case $k=0$ as a limiting case, even
though it never appears in the homotopy groups of spheres.
This phenomenon will recur many times below.
\end{rem}

\begin{rem}[{\bf Bocksteins}]\label{bockstein-1}Here is the another approach
to the calculation of Theorem \ref{1-line}.
There is a short exact sequence
of $G_1$-modules
\begin{equation}\label{ss1}
0 \to \ZZ/p(t) \mathop{\longr}^{\times p}
\ZZ/p^\infty(t) \longr \ZZ/p^\infty(t) \to 0
\end{equation}
which yields a long exact sequence in cohomology. There is
also split short exact sequence of groups (for $p > 2$)
$$
0 \to \ZZ_p \to G_1 \to \FF_p^\times \to 0
$$
and the action
of $G_1 = \ZZ_p^\times$ on $\ZZ/p(t) = \FF_p(t)$ is through
the powers of the evident action of $\FF^\times_p=C_{p-1}$ on
$\FF_p$. Thus
$$
H^\ast (G_1,\FF_p(\ast)) \cong \FF_p[v_1^{\pm 1}] \otimes H^\ast(\ZZ_p,\FF_p)
$$
where $\ZZ_p$ acts trivially on $\FF_p$. Thus
$$
H^\ast (G_1,\FF_p(\ast)) \cong \FF_p[v_1^{\pm 1}]  \otimes E(\zeta)
$$
where $\zeta \in H^1(G_1,\FF_p(0))$. In \cite{Rav} Theorem 6.3.21,
$v_1\zeta$ is called $h_{10}$.

To get the calculation
of $H^\ast(G_1,\ZZ/p^\infty(\ast))$ we need to calculate
the higher Bocksteins inherent in the long exact sequence
we have in cohomology. This computation requires essentially
the same information at the proof of Theorem \ref{1-line}.

At the prime $2$,  one can run the same argument, although
now the action of $G_1$ on $\FF_2(t)$ is trivial for all $t$ and
one has
$$
H^\ast(G_1,\FF_2(\ast)) \cong \FF_2[v_1^{\pm 1}] \otimes
\FF_2[\eta] \otimes E(\zeta)
$$
where $\eta$ and $\zeta$ are the evident generators of
$$
H^1(\{\pm 1\},\FF_2) \subseteq H^1(G_1,\FF_2(1))
$$
and
$$
H^1(\ZZ_2,\FF_2) \subseteq H^1(G_1,\FF_2(0))).
$$
These elements are variations of the elements $h_{10}$
and $\rho_1$ of \cite{Rav} Theorem 6.3.21.
\end{rem}

We now turn to the calculation of the $2$-line of the ANSS, at least for primes
$p > 3$. I'll get more sketchy than ever at this point, as the calculations are
beginning to get quite involved.

Because of Lemma \ref{0-line} and Theorem \ref{1-line} and the chromatic
spectral sequence we see that there is an isomorphism (in positive degrees)
$$
H^{2,\ast}(BP_\ast) \cong \mathrm{Ker}\{H^{0,\ast}(v_2^{-1}BP_\ast/I_2^\infty) \to 
H^{0,\ast}(v_3^{-1}BP_\ast/I_3^\infty)\}.
$$
Thus the first task is to compute
$$
H^{0,\ast}(v_2^{-1}BP_\ast/I_2^\infty) \cong H^0(G_2,R_2/\mm^\infty(t)).
$$
The module we wish to compute the cohomology of is very complicated (already here
at the $n=2$ case!). Some formulas can be found in \cite{DH}; however,
they aren't an immediate help. Thus we use the methods
of Bocksteins which I mentioned in Remark \ref{bockstein-1},
except we are now two Bocksteins away from the answer we
want. Rewrite
$$
R_2/\mm^\infty(t) = u^t R_2/\mm^\infty = u^tW(\FF_{p^2})[[u_1]]/(p^\infty,u_1^\infty).
$$
Then there are short exact sequences of $G_2$-modules
$$
0 \to u^t\FF_{p^2}[[u_1]]/(u_1^\infty)  \to u^tR_2/\mm^\infty
\mathop{\longr}^p u^tR_2/\mm^\infty \to 0
$$
and
$$
0 \to u^t\FF_{p^2} \mathop{\longr}^{v_1^{-1}}
u^{t-(p-1)}\FF_{p^2}[[u_1]]/(u_1^\infty)  \mathop{\longr}^{v_1}
u^t \FF_{p^2}[[u_1]]/(u_1^\infty)  \to 0.
$$
In the second of these, the element $v_1 = u^{p-1}u_1$ is
invariant under the $G_2$ action, so that we actually get $G_2$-module
homomorphisms. The action of $G_2$ on $u^t\FF_{p^2} = \FF_{p^2}(t)$
is through the quotient $\Gal(\FF_{p^2}/\FF_p) \rtimes \FF_{p^2}^\times$;
the cohomology we need was calculated in Corollary \ref{1cohG2}.
In particular we have identified two classes $h_0,h_1 \in
H^1(G_2,\FF_{p^2}(\ast))$.

With a little thought we see that $h_i$ is represented by
$t_1^{p^i} \in \Sigma(n)$. Because 
\begin{equation}\label{powert10}
t_1^{p^2} = v_2^{p-1}t_1
\end{equation}
(see Example \ref{Morava-stabilizer}) $h_i$, $i > 1$,
can rewritten as a power of $v_2$ times $h_0$ or $h_1$
and, in particular, is non-zero. For example
\begin{equation}\label{powert1}
h_2 = v_2^{p-1}h_0\qquad \mathrm{and}\qquad
h_3 = v_2^{p^2-p}h_1.
\end{equation}
These facts, and the
formula
\begin{equation}\label{etarv2}
\eta_R(v_2) = v_2 + v_1t_1^p-v_1^pt_1\qquad\mathrm{modulo}\qquad
(p)
\end{equation}
allows us to complete the Bockstein calculation to decide
how divisible the class $v_2^k$ becomes in the
cohomology of $\FF_{p^n}[[u_1]]/(u_1^\infty)$. One way
to phrase the following result is that $1 = v_2^0$ is
infinitely $v_1$ divisible and if $k = p^ik_0$ with
$(p,k_0) = 1$, then $v_2^k$ is $v_1^{a_i}$ divisible
where $a_i \to \infty$ as $i \to \infty$.
This is analogous to the $p$-divisibility 
in Theorem \ref{1-line}: the case $k=0$ is the limiting case
achieved as $k \to 0$ in the $p$-adic topology. See Remark
\ref{p-adic-nghbd}.
Note that since $p=0$ in $\FF_{p^n}[[u_1]]/(u_1^\infty)$,
the cohomology groups of this module are modules
over $\FF_p[v_1]$.

In the following result, the element $x_i^{s}/v_1^j$ is the corrected
version of the element $v_2^{p^is}/v_1^j$ needed to make the Bocksteins
work out. See Remark \ref{where-bock}.

\begin{thm}\label{2-line-cont}In the $\FF_2[v_1]$-module
$$
H^0(G_2,\FF_{p^n}[[u_1]]/(u_1^\infty)(\ast))
$$
is spanned by the following sets of non-zero elements, with the 
evident $v_1$ multiplication:
\begin{enumerate}

\item $1/v_1^j$, $j \geq 1$;

\item if $s = p^is_0$ with $(p,s_0)=1$, there is an element $x_i \in E(2)_\ast$
so that $x_i \equiv v_2^{p^i}$ modulo $(p,v_1)$ and we have the
classes 
$$
x_i^{s_0}/v_1^j
$$
where
$$
1 \leq j \leq a_i= \brackets{1,}{i=0;}{p^i+p^{i-1} -1,}{i \geq 1.}
$$
\end{enumerate}
\end{thm}

\begin{rem}\label{where-bock}
This is Theorem 5.3 of \cite{MRW}.
I won't give a proof, but it's worth pointing
out how the number $a_i$ arises. Write $d^n$ for the $n$th $v_1$-Bockstein.
Then the formula \ref{etarv2} immediately implies that
$$
d^1(v_2^s) = sv_2^{s-1}h_1.
$$
Thus, if $s$ is divisible by $p$, $d_1(v_2^s) = 0$ and we compute
also using $\ref{etarv2}$ and \ref{powert1} that the next Bockstein is
$$
d^p(v_2^s) = (p/s) v_2^{p(s-1)}h_2 = (p/s) v_2^{s-1}h_0.
$$
Thus, if $p^2$ divides $s$ we continue and we have, now
again using \ref{powert1}
\begin{align*}
d^{p^2}(v_2^s) &= (s/p^2)v_2^{s-p^2}h_3\\
&= (s/p^2)v_2^{s-p}h_1= d^1((s/p^2)v_2^{s-p+1}).
\end{align*}
Thus $d^{p^2}(v_2^s) = 0$ in the Bockstein spectral sequence and,
making the correction dictated by this equation we calculate that
$d^{p^2+p-1}(v_2^s) \ne 0$.
\end{rem}

The next step is to calculate the groups
$$
H^0(G_2,R_2/\mm^\infty(\ast))
$$
using the the Bockstein based on multiplication by $p$. Put another
way, we must decide how $p$-divisible the elements of Theorem
\ref{2-line-cont} become. Since this
requires knowledge of $H^1(G_2,\FF_{p^2}[[u_1]]/(u_1^\infty)(\ast))$,
I won't even pretend to outline the calculation. However, here
are two comments.
\begin{enumerate}

\item The first differential in the chromatic spectral sequence
$$
H^{0,\ast}(v_1^{-1}BP_\ast/(p^\infty)) \to H^{0,\ast}(v_2^{-1}BP_\ast/I_2^\infty)
$$
must be an injection in negative degrees, this forces a maximum $p$ divisibility
on the elements $1/v_1^j$ and we get a ``reverse'' image of $J$ pattern
working downward from $1/v_1$.

\item The case when $s=p^is_0$ -- that is, of the elements $x_i^{s_0}/v_1^j$ --
must approach the case of $1/v_1^j$ as $i \to \infty$. Thus we would
expect the beginning of a reverse image of $J$ pattern working downward
from $x_i^{s_0}/v_1$. This pattern is not complete, however, by which
I mean we get the image of $J$ divisibility if $j$ is small, but as
$j$ becomes large, the divisibility becomes constrained. 
\end{enumerate}

One reason for the constraints in (2) was suggested to me by Mark Mahowald.
In the clASS, multiplication by $v_1$ and multiplication by $p$ raise filtration.
Thus, after dividing
by $v_1$ a certain of times, we have arrived at a low filtration in the clASS. 
Since we expect image of $J$ type divisibility to occur in
relatively high filtration in the clASS, we don't
expect unconstrained divisibility by $p$ at
this point. The exact result is Theorem 6.1 of \cite{MRW}:

\begin{thm}\label{2-line-cont-2}Let $p \geq 3$. The groups 
$$
H^0(G_2,R_2/\mm^\infty(\ast))
$$
are spanned by the non-zero classes, with evident multiplication by $p$,
\begin{enumerate}

\item $1/p^{k+1}v_1^j$ where $j > 0$ and $p^k$ divides $j$; and

\item if $s = p^is_0$ with $(p,s_0)=1$, we have the elements
$$
x_i^{s_0}/p^{k+1}v_1^j
$$
with $1 \leq j \leq a_i$ and with $k$ subject to the twin requirements that
$$
p^k | j\qquad\mathrm{and}\qquad a_{i-k-1} < j \leq a_{i-k}.
$$
\end{enumerate}
\end{thm}

The first of these contraints is the image of $J$ divisibility; when
$j$ is small, the second contraint won't apply. On the other hand,
if $j$ is large, the second constraint dominates. For example, if
$$
p^{n-1} + p^{n-2} - 1 = a_{n-1} < j \leq a_n = p^n + p^{n-1} -1
$$
then we only have elements of order $p$.

With the result in hand, the question then remains how many of the
elements are non-zero permanent cycles in the chromatic spectral
sequence. None of the elements in negative degree can be; thus
the limiting case of the reverse image of $J$ pattern cannot
survive, even though it has avatars throughout the rest of the groups. 
There is another class
of elements that cannot survive as well; here is heuristic argument
for this case.

Consider the family of elements
$$
x_i/pv_1^j,\qquad i \geq 2.
$$
These are roughly $v_2^{p^i}/pv_1^j$. The element $x_i/pv_1^{p^i}$
does survive to an element
$$
\beta_{p^i/p^i} \in \Ext_{BP_\ast BP}^2(\Sigma^{2p^i(p-1)}BP_\ast,BP_\ast)
$$
which reduces to a non-zero class
\begin{equation}\label{odd-arf}
b_i \in \Ext^2_{A_\ast}(\Sigma^{2p^i(p-1)}\FF_p,\FF_p)
\end{equation}
in the clASS -- the odd primary analog of the Kervaire invariant
class.\footnote{We would expect such elements to be primordial 
-- and certainly not $v_1$-divisible.} On the general principal,
alluded to above,  that
dividing by $v_1$ cannot
raise the cohomological degree $s$, we see immediately that
we've run out of room --  there's no place to put $x_i/pv_1^j$, $j > p^i$.
The following result follows from Lemma 7.2 of \cite{MRW}.

\begin{thm}\label{2-line}Suppose $s = p^is_0 > 0$ with $(p,s_0) = 1$.
Then all of the classes $x_i^{s_0}/p^{k+1}v_1^j$ survive to
$E_\infty$ in the chromatic spectral sequence {\bf except}
$$
x_i/pv_1^j,\qquad 2 \leq i ,\ p^i < j \leq a_i.
$$
If $s < 0$, the elements $x_i^{s_0}/p^{k+1}v_1^j$ do not
survive.
\end{thm}

We write
\begin{equation}\label{beta-names}
\beta_{s/k,j} \in \Ext^2_{BP_\ast BP}(\Sigma^\ast BP_\ast,BP_\ast)
\end{equation}
for the unique class detected by $x_i^{s_0}/p^{k+1}v_1^j$. If
$k = 0$ we write $\beta_{s/j}$; if $j=1$ as well, then we recover
$\beta_s$ of Theorem \ref{period1}.

\begin{rem}\label{survivors} Now, of course, one can ask which
of these extended $\beta$-elements are infinite cycles in the 
ANSS. A survey of known results is in \cite{Rav} \S 5.5. For
example, it is a result of Oka's that enough of the
generalized $\beta$s do detect homotopy class to show
that for all $n$ there is a $k$ so that $\pi_kS^0$ contains
a subgroup isomorphic to $(\ZZ/p)^n$. 

There are also some spectacular non-survival results. For example,
Toda's differential reads for primes $p \geq 3$ 
$$
d_{2p-1}\beta_{p/p} = \alpha_1\beta_1^{p-1}.
$$
up to non-zero unit. Ravenel used an induction argument to show
$$
d_{2p-1} \beta_{p^{i+1}/p^{i+1}} = \alpha_1\beta_{p^i/p^i}^{p-1}
\qquad \mathrm{modulo}\qquad \mathrm{Ker}\ \beta_1^{c_i}
$$
up to non-zero unit. Here $c_i = p(p^i-1)/(p-1)$. This result
has as a corollary that if $p > 3$, the elements
$b_i$ of \ref{odd-arf} cannot be permanent cycles in the
clASS. At the prime $3$ the divided $\beta$s can gang
up to produce infinite cycles that map to $b_i$ -- some
some $b_i$ survive and some do not. Assuming the
Kervaire invariant elements at the prime $2$ all do 
exist, this is an example of the Mahowald axiom that
``something has to be sick before it dies''.

That the
differential is as stated is fairly formal consequence of Toda's
differential; the subtle part is
to show that the target is non-zero modulo the listed 
indeterminacy. To do this, Ravenel uses the
chromatic spectral sequence and detects the necessary elements
in
$$
H^\ast (G_{n},\FF_{p^n}(\ast)).
$$
with $n=p-1$.
In this case, the Morava stabilizer group $S_{p-1} \subseteq G_{p-1}
=G_n$
has an element of order $p$; choosing such an element gives a map
$$
H^\ast (G_{n},\FF_{p^{n}}(\ast)) \to H^\ast(\ZZ/p,\FF_{p^n}(\ast))
$$
and the necessary element is detected in the (highly-non-trivial)
target. A variation on this technique was used in \cite{nave} to get
non-existence results for the Smith-Toda complexes $V(p-1)$.
\end{rem}

\end{document}